\newif\ifarXiv         
\newif\ifjournal        
\journalfalse  \arXivtrue       
\documentclass[11pt]{article}
\usepackage[small,compact]{titlesec}
\usepackage{amsmath, amssymb, amsfonts, mathabx, setspace}
\usepackage{graphicx,epsfig,wrapfig,caption,epsfig,subcaption,sidecap}
\usepackage{url,color, verbatim,algorithmicx}
\usepackage[noend]{algpseudocode}
\usepackage[ruled,boxed]{algorithm} 

\usepackage{enumitem} 
\usepackage{booktabs}  
\usepackage[sort,compress]{cite}
\usepackage[scaled]{helvet}
\usepackage[T1]{fontenc}
\usepackage[bookmarks=true, bookmarksnumbered=true, colorlinks=true,   pdfstartview=FitV,
linkcolor=blue, citecolor=blue, urlcolor=blue]{hyperref}
\usepackage[english]{babel}  
\usepackage{multirow}
\let\OLDthebibliography\thebibliography
\renewcommand\thebibliography[1]{
  \OLDthebibliography{#1}
  \setlength{\parskip}{0pt}
  \setlength{\itemsep}{3pt}
}

\usepackage{fullpage}

\def\pitilde{\tilde{\pi}}

\def\calA{\mathcal{A}}

\def\bA{\mathbf{A}}    
\def\bb{\mathbf{b}}

\def\bx{\mathbf{x}}
\def\bB{\mathbf{B}}

\def\bT{\mathbf{T}}
\def\bP{\mathbf{P}}

\def\R{\mathbb{R}}                               
\def\P{\mathbb{P}}    \def\E{\mathbb{E}}

\newcommand{\prob}[1]{\mathbb{P}\left\{{#1}\right\} }

\newcommand{\floor}[1]{\lfloor{#1}\rfloor}
\newcommand{\argmin}[1]{\underset{#1}{\operatorname{arg}\operatorname{min}}\;}

\newtheorem{theorem}{Theorem}

\newtheorem{assumption}[theorem]{Assumption}

\newtheorem{corollary}[theorem]{Corollary}
\newtheorem{definition}[theorem]{Definition}
\newtheorem{example}[theorem]{Example}
\newtheorem{lemma}[theorem]{Lemma}
\newtheorem{proposition}[theorem]{Proposition}
\newtheorem{remark}[theorem]{Remark}

\newenvironment{proof}[1][Proof]{\noindent\textbf{#1.} }{\ \rule{0.5em}{0.5em}}

\numberwithin{theorem}{section}
\numberwithin{equation}{section}

\usepackage{soul}

\definecolor{silver}{RGB}{192,192,192}

\usepackage{tikz}                 

\definecolor{state1}{RGB}{255,255,0}
\definecolor{state2}{RGB}{0,206,209} 
\definecolor{state3}{RGB}{0,0,139} 

\usepackage{authblk}

\title{
Probabilistic cellular automata with local transition matrices: synchronization, ergodicity, and inference 
}

\author[a]{Erhan Bayraktar \thanks{erhan@umich.edu}}
\author[b]{Fei Lu \thanks{feilu@math.jhu.edu}} 
\author[b,c]{Mauro Maggioni \thanks{mauromaggionijhu@icloud.com}}
\author[d]{Ruoyu Wu \thanks{ruoyu@iastate.edu}}
\author[c]{Sichen Yang \thanks{syang114@jhu.edu}}
 
\affil[a]{Department of Mathematics, University of Michigan, Ann Arbor, USA.}
\affil[b]{Department of Mathematics, Johns Hopkins University, Baltimore, USA.} 
\affil[c]{Department of Applied Mathematics and Statistics, Johns Hopkins University, Baltimore, USA.}
\affil[d]{Department of Mathematics, Iowa State University, Ames, USA.}

\date{} 
  
\begin{document}

\vspace{-6mm}
\maketitle
\vspace{-10mm}
\begin{abstract}
We introduce a class of probabilistic cellular automata that are capable of exhibiting rich dynamics such as synchronization and ergodicity and can be easily inferred from data. The system is a finite-state locally interacting Markov chain on a circular graph. Each site's subsequent state is random, with a distribution determined by its neighborhood's empirical distribution multiplied by a local transition matrix. We establish sufficient and necessary conditions on the local transition matrix for synchronization and ergodicity. Also, we introduce novel least squares estimators for inferring the local transition matrix from various types of data, which may consist of either multiple trajectories, a long trajectory, or ensemble sequences without trajectory information. Under suitable identifiability conditions, we show the asymptotic normality of these estimators and provide non-asymptotic bounds for their accuracy. 
\end{abstract}



\setcounter{tocdepth}{1}
\tableofcontents

\section{Introduction}
Interacting systems, including probabilistic cellular automata (PCA) \cite{toom1994critical,lebowitz1990statistical,louis2018probabilistic} and interacting particle systems (IPS) \cite{liggett1985interacting,durrett2007random,aldous2013interacting,grimmett2018probability}, have a wide range of applications in Physics, Computer Science, Electrical Engineering, Economics, Biology, among others. These applications have driven a growing interest in studying the dynamics of these systems and inferring model parameters from observational data. As Aldous \cite{aldous2013interacting} pointed out, it is ``the most broad-ranging currently active field of applied probability''; however, ``it is easy to invent and simulate models, but hard to give rigorous proofs or to relate convincingly to real-world data''. 

In this work, we introduce a class of PCA that exhibits rich dynamics, yet can be easily inferred from observational data. We rigorously prove dynamical properties such as synchronization and ergodicity, and then construct computationally efficient estimators of the model parameters, for which we prove desirable properties such as asymptotic normality and non-asymptotic bounds.         

\begin{figure}
        \centering
        \includegraphics[width=0.27\textwidth]{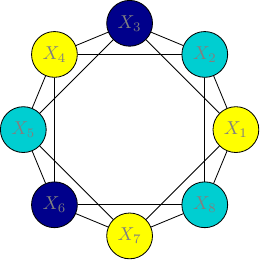} 
\includegraphics[width=0.34\textwidth]{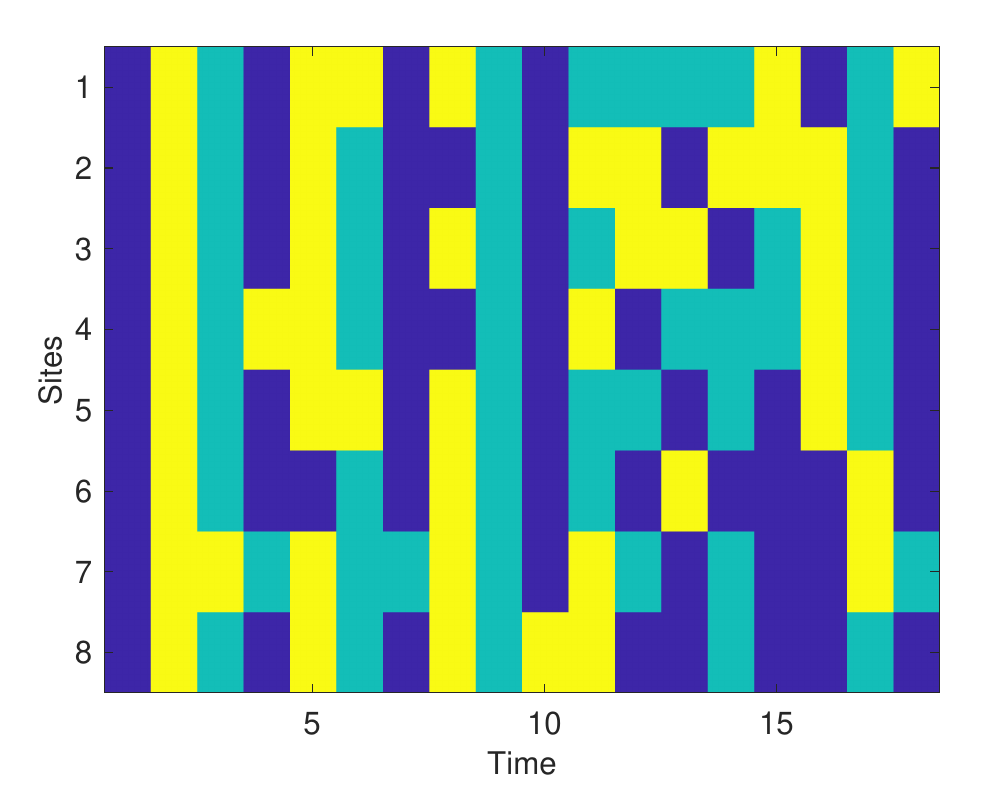} 
        \includegraphics[width=0.34\textwidth]{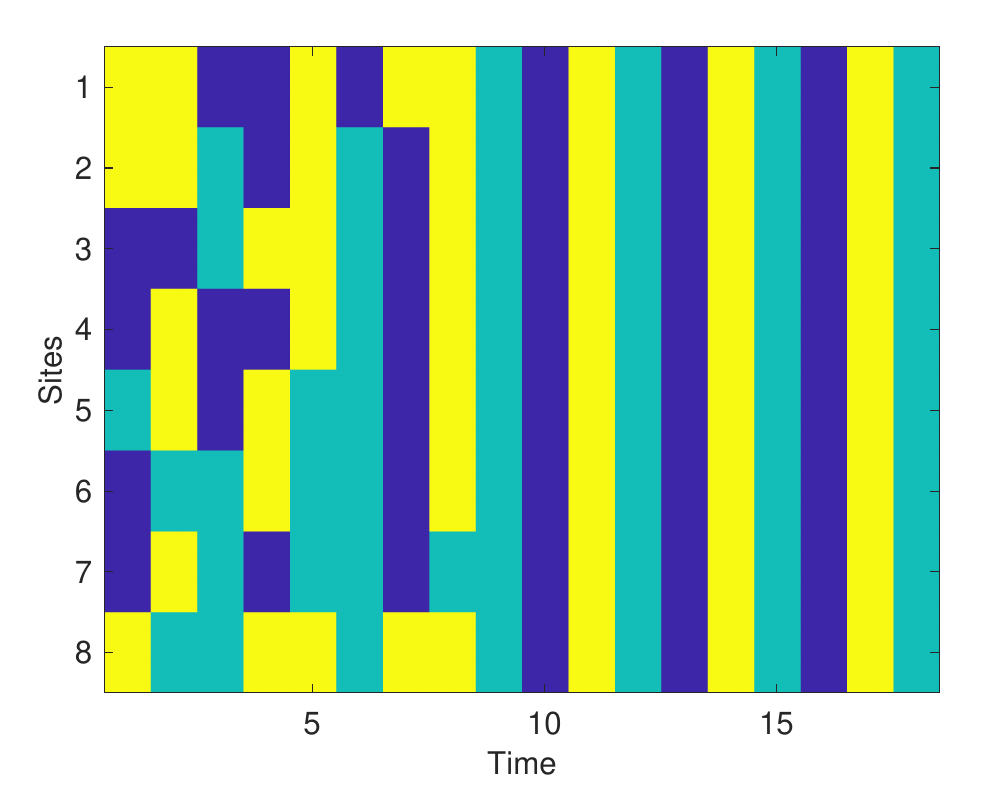} \\
        (a) PCA on a circle \hspace{1.5cm} (b) From deterministic to stochastic  \hspace{1.5cm}  (c) Synchronization\\
\caption{(a) The system of $N= 8$ agents on a graph with an alphabet $\calA= \{\text{Yellow, Turquoise, and Blue}\}$ with $K=3$. Each agent's transition depends linearly on the empirical distribution of its nearest neighbors with $n_v=2$ agents on each side. (b) The system moves from deterministic to stochastic dynamics; see Example \ref{exp:T_move2next}. (c) The system transitions from stochastic to deterministic dynamics, achieving a synchronization; see Example \ref{exp:T-permutation}. }\label{fig:circleGraph}
\end{figure}

This class of PCA is defined on an $N$-node cyclic graph $(V, E)$, with a finite alphabet $\calA$, and every site updates independently with a distribution determined by the empirical distribution of its neighborhood multiplied by a local transition matrix.  
Each PCA is a finite-state Markov chain
\[
X(t) =(X_1(t),\ldots,X_N(t))=: X_{1:N}(t) \in \calA^{N}, \quad \calA = \{1,\ldots, K\} =: [K] 
\] 
on a cyclic graph $(V,E)$ with nodes indexed by $V= \{1,2,\dots,N\} =: [N]$ and with edges in $E$ connecting nodes within distance $n_v$, that is, nodes $n$ and $n'$ are connected whenever $|n-n'| \le n_v$ (modulo $N$); see Figure \ref{fig:circleGraph}(a) for an illustration.  
Conditional on $X(t)$, each vertex $n$ makes updates independently depending on its neighborhood 
$$V_n= \{n-n_v,\ldots, n, n+1,\ldots, n+n_v\}\quad \text{(modulo $N$)}\,,$$ 
and, since the state updates at each vertex are simultaneous and independent, the transition probability of $X(t)$ is in the form  
\begin{equation}\label{eq:MC}
\prob{ X(t+1)=(x_1,\ldots,x_N) | X(t) } = \prod_{n= 1}^N\prob{X_n(t+1) = x_n| (X_{n'}(t))_{n'\in V_n} }.
\end{equation}
In our model we assume that the local transition probability $\prob{X_n(t+1) = x_n| (X_{n'}(t))_{n'\in V_n} }$ of the $n$-th vertex depends \emph{linearly} on the empirical distribution of its neighborhood:
\begin{equation}\label{eq:PT}
\prob{ X_n(t+1) = x_n | (X_{n'}(t))_{n'\in V_n}  } 
 =    \varphi_{n,t}\bT(\cdot,x_n),   
\end{equation}
where the \emph{local transition matrix} $\bT\in [0,1]^{K\times K}$ is row-stochastic (i.e., $\sum_{k=1}^K \bT(j,k) = 1$, $j\in[K]$),
and $\varphi_{n,t}\in \R^{1\times K}$ is the \emph{local empirical distribution} of the vertex's neighborhood $V_n$ at time $t$:
 \begin{equation}\label{eq:phi}
\varphi_{n,t} = (\varphi_{n,t}(1),\ldots, \varphi_{n,t}(K)), \quad \text{ with }   
\varphi_{n,t}(k) := \frac{1}{|V_n|} \sum_{i\in V_n} \delta_{X_i(t)}(k),  \quad k\in[K]\,, 
\end{equation} 
where $|V_n|$ is the cardinality of $V_n$, and $\delta$ is the Kronecker delta function. 
Equivalently, the Markov chain $X(\cdot)$ has a \textit{global transition matrix} $\bP\in \R^{K^N\times K^N}$ determined by $\bT$ by
\begin{equation}\label{eq:T2P}
\bP(x_{1:N},y_{1:N}) = \prod_{n=1}^N  \varphi_n\bT_{\cdot,y_n} = \prod_{n=1}^N  \frac{1}{|V_n|}\sum_{i\in V_n} \bT_{x_i,y_n} 
\end{equation}
for any configurations $x_{1:N}, y_{1:N}\in [K]^N$, where $\varphi_n(k):= \frac{1}{|V_n|}\sum_{i\in V_n} \delta_{x_i}(k)$ for $k\in [K]$ and $n\in [N]$. 
The evolution of the system may also be described as follows: at each update step, every vertex samples a state from its neighbors uniformly at random, and then independently jumps according to the local transition matrix $\bT$. 

The key feature of our model is the {\em{linear}} dependence of the local transition probability in \eqref{eq:PT} on the local empirical distribution, represented by a local transition matrix. Such a linear dependence significantly reduces the number of parameters describing the local transition probability, which has a size $K^{|V_n|+1}$ since it assigns an $ \R^{1\times K}$-valued probability to each of the $K^{|V_n|}$ possible states of the neighborhood $V_n$. Without the linear dependence, the transition probability is overly complicated for analysis and requires a significant amount of data for its estimation. In contrast, with the linear dependence, our model has only $K^2$ parameters instead of $K^{|V_n|+1}$, significantly reducing the model complexity and the amount of data needed for inference.

Yet, the Markov chains in these PCA can exhibit rich dynamics such as synchronization and ergodicity. Figure \ref{fig:circleGraph} illustrates two systems with $(N, K,n_v) = (8,3,2)$ and with different local transition matrices. Figure \ref{fig:circleGraph}(b) shows a trajectory exhibiting a transition from deterministic to stochastic dynamics, and Figure \ref{fig:circleGraph}(c) shows a trajectory exhibiting synchronization. 

\subsection{Main results}
\paragraph{Synchronization and ergodicity.} 
We first study how the local transition matrix $\bT$ determines important features of the global dynamics by establishing sufficient and necessary conditions on $\bT$ for synchronization and ergodicity. We have a full characterization of the dynamics when $\bT$ is irreducible: (i) Proposition \ref{prop:periodicity} shows that if $\bT$ had period $d<K$, then a state of $X(t)$ either has a period $d$ or is transient; 
(ii) Theorem \ref{thm:sync} shows that the system will achieve a synchronization if and only if $\bT$ is periodic with period $K$;  and (iii) Theorem \ref{thm:ergodic}
shows that the system is exponentially ergodic if and only if $\bT$ is aperiodic. 
In short, 
\[
\bT \text{ irreducible}: 
\begin{cases}
 \text{ periodic with period $d<K$} &\Leftrightarrow \text{a state of $X$ is periodic or transient}  \\
\text{ periodic with period $K$} &\Leftrightarrow \text{synchronization} \\
\text{ aperiodic} &\Leftrightarrow \text{exponential ergodicity}.   	
 \end{cases}
\]

\paragraph{Global and local transition matrices.} We then study the dependencies between the global transition probability matrix $\bP\in \R^{K^N\times K^N}$ of the Markov chain and the local transition matrix $\bT\in\R^{K\times K}$. Theorem \ref{thm:1-1PT} shows that there is a 1-1 map between them; additionally, the invariant measure of $\bT$ is the marginal distribution of the invariant measure of $\bP$.  Theorem \ref{thm:Lip_P_T} shows that $\bP$ is Lipschitz in $\bT$ under the total variation norm. 

\paragraph{Inference of $\bT$.} In Section \ref{sec:inference} we introduce least squares estimators (LSEs) for inferring the local transition matrix from data, and study their properties. We consider different types of observational data: it may consist of multiple trajectories, a single long trajectory, or ensemble sequences without trajectory information. Except for the case of a single long trajectory, the system may be non-ergodic. These LSEs use the marginal distributions of each vertex and are more efficient to compute than the maximal likelihood estimator, which would involve non-convex optimization.  We specify \emph{identifiability} conditions for these LSEs and a non-identifiability for inference using the stationary distribution in Section \ref{sec:infer_from_inv}. We show the \emph{asymptotic normality} of these estimators in Theorems \ref{thm:AN_lse} and \ref{thm:AN_ensLSE}; we provide \emph{non-asymptotic bounds} for their accuracy in Section \ref{sec:non-asymp}.  Numerical tests show that the LSE with trajectory information is more accurate than the LSE without trajectory information, while both converge at the optimal rate $M^{-1/2}$, with $M$ being the sample size.

\subsection{Related work}

\paragraph{Probabilistic cellular automata (PCA).} PCA are a large class of interacting discrete stochastic dynamical systems for the modeling of a wide range of physical and societal phenomena; we refer the reader to \cite{moran1958random,toom1994critical,lebowitz1990statistical,louis2018probabilistic} and references therein for applications. There is consistent interest in studying the ergodicity of the systems; see, e.g., Dawson \cite{dawson1975information} for a system with interacting subsystems, Follmer and Horst \cite{follmer2001convergence} for the averaged process of an interacting Markov chain with an infinite set of sites, B\'erard  \cite{berard2023coupling} for the exponential ergodicity of a 1D PCA with a local transition kernel on a three-state alphabet, and Casse \cite{casse2023ergodicity} for the ergodicity of a PCA with binary alphabet via random walks. The innovation in our model above is introducing a local transition matrix, which enables efficient estimation, while leaving the system capable of producing rich dynamics, exhibiting phenomena such as synchronization or ergodicity.

\paragraph{Interacting particle systems (IPSs).} IPSs are continuous-time Markov processes on certain spaces of configurations of finitely or infinitely many interacting particles. The state space can be either discrete, such as the stochastic Ising model or the voter model   
\cite{liggett1985interacting,durrett2007random,aldous2013interacting,grimmett2018probability}, or continuous in the form of stochastic differential equations \cite{cucker2008_FlockingNoisy,cattiaux2018_StochasticCucker,lacker2021locally}. The interaction rules, either short-range or long-range, are often specified by functions called interaction kernels/potentials \cite{liggett1985interacting,cattiaux2018_StochasticCucker} or rate functions \cite{aldous2013interacting}. Thus, our local transition matrix can be viewed as a counterpart of these interaction kernels or rate functions.    

\paragraph{Inference of the local transition matrix.} The inference of the local transition matrix is akin to the nonparametric estimation of the interaction kernel of interacting particle systems in \cite{LZTM19pnas,LMT21_focm,della2022nonparametric,LWLM24}, where inference leads to a linear inverse problem and is solved by least squares. However, the estimators in those works maximize the likelihood; here, our LSEs are different from the maximal likelihood estimator, which would lead to a constrained non-convex optimization problem, as discussed in Section \ref{sec:lse_traj}. Also, while the identifiability conditions are specified based on the large sample limit case, in the same spirit as the coercivity condition on function spaces in \cite{LLMTZ21,LiLu23coercivity}, this study considers parametric inference, so the identifiability conditions are less restrictive.

We use the notations in Table \ref{tab:notation} throughout the paper. We denote the entries of $\bT\in \R^{K\times K}$ by $\bT_{jk}$ with $j,k\in [K]$, and denote the entries of $\bP\in \R^{K^N\times K^N}$ by $\bP(x,y)$ with $x,y\in [K]^N$, where $x=x_{1:N}:= (x_1,\ldots, x_N)$ with $x_i\in [K]$.  
\begin{table}
\caption{Notations. 
} \label{tab:notation} 
\begin{tabular}{l |  l  } 
\hline
$ k\in \calA=[K]= \{1,\ldots, K\}$ & alphabet set for state values  \\ 
$n\in [N]= \{1,2,\ldots, N\}$ & index of vertices/agents in the graph \\
$m\in [M], \ t\in [L]$ & index of sample trajectories and index of time \\ \hline 
$X(t) \in \calA^N$ & state of the Markov chain at time $t$\\
$V_n$ and $n_v$    & vertex $n$'s neighborhood $V_n$, consisting of $2n_v+1$ vertices  \\
$ \varphi_{n,t} =(\varphi_{n,t}(1), \ldots, \varphi_{n,t}(K))$ & empirical distribution in $V_n$ at time $t$; $[0,1]^{1\times K}$-valued \\
$c_{n,t} = (c_{n,t}(1),\ldots,c_{n,t}(K))  $ & empirical distribution of $X_n(t)$: $ c_{n,t}(k) := \delta_{X_n(t)}(k)$ 
  \\ \hline
$\bT\in [0,1]^{K\times K}$ &  local transition matrix: $\sum_{k=1}^K \bT_{jk} = 1, \forall j \in [K]$  \\
$\bP \in [0,1]^{K^N \times K^N}$ & (global) transition matrix  $ \big(\bP(x,y)\big)_{x,y\in \calA^N}$\\
$\|\cdot\|$, $\|\cdot\|_F$, $\|\cdot\|_{op}$, $\|\cdot\|_p$ & Euclidean,  Frobenius, operator, and $\ell_p$ norms \\
\hline
\end{tabular}
\end{table}


\section{Dynamical properties: synchronization and ergodicity} 
\label{s:example:lattice}
This section studies the dynamical properties of the process $X(\cdot)$ in \eqref{eq:MC}-\eqref{eq:PT} as a Markov chain with $K^N$ states. We characterize the long-time behavior of the system when the local transition matrix $\bT$ is irreducible: Theorem {\rm\ref{thm:sync}} shows that the system achieves a synchronization if and only if $\bT$ is periodic with period $K$; Proposition {\rm\ref{prop:periodicity}} shows that the system will eventually be periodic with the same period as that of $\bT$; and Theorem {\rm\ref{thm:ergodic}} and Proposition {\rm \ref{prop:irreduciblePT}} show that the system is exponentially ergodic if and only if $\bT$ is aperiodic.

We recall the following basic notions about a finite-state Markov chain, denoted by $Z(t)$, with states $[n]$ and transition matrix $A\in \R^{n\times n}$. 
\begin{itemize}
\item The transition matrix $A\in \R^{n\times n}$ is called \emph{irreducible} if $\forall i,j\in [n]^2$, $\exists t \in \mathbb{N}$, such that $(A^t)_{ij}>0$. 

\item The \emph{period} $d(i)$ of a state $i$ is the greatest common divisor of all $m$ such that $(A^m)_{ii} > 0$, i.e., $d(i) := gcd\{ m\in \mathbb{N}: (A^m)_{ii} >0\}$. 
When $A$ is irreducible, the period is the same for all states, and it is called the period of $A$. 
An irreducible matrix $A$ of period $1$ is called \emph{aperiodic}.

\item A state $i$ is \emph{recurrent} if $\prob{\tau_i < \infty | Z(0)=i} = 1$, where $\tau_i = \min\{ t \in \mathbb{N}: Z(t) = i \}$; in other words, the chain started at $i$ returns to $i$ in finite time with probability one.  State $i$ is \emph{transient} if $\prob{\tau_i < \infty | Z(0)=i} < 1$. 
State $i$ is positive recurrent if $\mathbb{E}[\tau_i | Z(0)=i] < \infty$.
\end{itemize}

In the following, we first present a few examples of PCA models that exhibit a variety of dynamical properties. Then, we study the sufficient and necessary conditions for the system to synchronize or to be ergodic.

\subsection{Examples: stochastic dynamics and synchronization}
We introduce four examples: non-interacting agents, the smallest model, a system transitioning from deterministic to stochastic dynamics, and a system achieving synchronization. 

\begin{example}[Non-interacting agents] When the agents do not interact, i.e., $n_v=0$, they move independently according to a Markov chain with $\bT$ as the probability transition matrix. That is, the process $X(t) = (X_1(t),X_2(t),\cdots, X_N(t))$ is a vector of $N$ independent Markov chains, all governed by the same transition matrix $\bT$. 
\end{example}

\begin{example}[Smallest model: $(N,K)=(2,2)$] 
\label{exp:full-nbhd-N2K2}
Consider the model with $(N,K,n_v)=(2,2,1)$. Since the neighborhood is the full network for each site, the local empirical distributions are the same for all sites, that is, we have $\varphi_1 = \varphi_2$ for all states. The following table shows the local empirical distributions and the global transition matrix with a local transition matrix  $\bT= \begin{pmatrix} \bT_{11}& \bT_{12} \\ \bT_{21}& \bT_{22} \end{pmatrix}$: 
  \[ \begin{array}{c|cc | cccc}\hline
   \, &  &  & \multicolumn{4}{c}{\text{ Global transition probability}\, \bP } \\  \hline
  x  & & \varphi_1 = \varphi_2 & 11   & 12  &  21  & 22 \\ \hline
 11 &  & (1, 0) &\bT_{11}^2 & \bT_{11}\bT_{12} &  \bT_{11}\bT_{12} & \bT_{12}^2 \\
 12 &  & (0.5, 0.5)&\quad  \frac{1}{4} (\bT_{11}+\bT_{21})^2 \quad  & \quad  \bP(12,12)  \quad  & \quad \bP(12,21) \quad  &  \frac{1}{4} (\bT_{12}+\bT_{22})^2  \\
 21 &  & (0.5, 0.5)&   \frac{1}{4} (\bT_{11}+\bT_{21})^2 &  \bP(21,12)  &  \bP(21,21)  &  \frac{1}{4} (\bT_{12}+\bT_{22})^2 \\
 22 &  &  (0, 1) & \bT_{21}^2 & \bT_{11}\bT_{12} &  \bT_{11}\bT_{12} & \bT_{22}^2 \\ \hline 
 \end{array} \]
 where $ \bP(12,12)=\bP(12,21) = \bP(21,12) = \bP(21,21)= \frac{1}{4} (\bT_{11}+\bT_{21}) (\bT_{12}+\bT_{22}) $. 
\end{example}

  \begin{example}[From deterministic to stochastic dynamics]\label{exp:T_move2next} 
  The system can change from deterministic to stochastic dynamics.  Let $\bT$ have zero entries except $\bT_{k,k+1}=1$ for $k=1,\dotsc,K-1$ and $\bT_{K,k}=\frac{1}{K}$ for $k\in [K]$; that is, 
 \begin{equation*}
 	 \bT= \begin{bmatrix} 
   0 & 1 & 0 & \cdots & 0\\
   0 & 0 & 1 & \cdots &0 \\
   \vdots & \vdots & \ddots &\ddots &\vdots\\
   0   & 0  &0 & \cdots & 1 \\
   1/K &   1/K &  1/K & \cdots &  1/K
\end{bmatrix} 
.\end{equation*}
The system's dynamics will move from deterministic to stochastic if it starts with state $(1,\dotsc,1)\in \calA^N$. Specifically, note that we have $\varphi_{n,0}= (1,0,\cdots,0)=: e_1\in \R^{1\times K}$ for each $n\in [N]$. Then, the value of each site moves to the next value deterministically, i.e., $X(t)= (t+1,t+1,\cdots, t+1)\in \calA^{N}$, for $t\leq K-1$. Correspondingly, we have $\varphi_{n,t}= e_{t+1}\in \R^{1\times K}$ for each $n\in [N]$ and $t\leq K-1$.  When $t>K-1$, the move becomes stochastic as the last row of the local transition matrix injects randomness. Figure {\rm\ref{fig:circleGraph}(b)} shows a typical trajectory of the system with $(N,K,n_v)=(8,3,2)$. 
 \end{example}

\begin{example}[Synchronization: from stochastic to deterministic dynamics] \label{exp:T-permutation}
When $\bT$ is a permutation matrix such that it is irreducible with period $K$, e.g., 
\begin{equation*}
 \bT= \begin{bmatrix} 
   0 & 1 & 0 & \cdots & 0\\
   0 & 0 & 1 & \cdots &0 \\
   \vdots & \vdots & \ddots &\ddots &\vdots\\
   0   & 0  &0 & \cdots & 1 \\
   1 &   0 & 0 & \cdots &  0
\end{bmatrix}
, \end{equation*}
the Markov chain will achieve synchronization (see Theorem {\rm\ref{thm:sync}}), in which all sites move from one state to another with the same period as $\bT$, as demonstrated in Figure {\rm\ref{fig:circleGraph}(c)} for a typical trajectory of the system with $(N,K,n_v)= (8,3,2)$. In particular, the deterministic dynamics of the Markov chain after synchronization is as follows. Without loss of generality, suppose that it starts from the state $X(0)= (1,\dotsc,1)\in \calA^N$. The local empirical distributions are $\varphi_{n,0}= (1,0,\cdots,0)=: e_1\in \R^{1\times K}$ for each $n\in [N]$. Then, all vertices move uniformly from one state to the next, that is, $X(t)= (t+1,\dotsc,t+1)$ for $t \leq K- 1$, and then repeat periodically, as shown in the following tabular. 
 \[
\begin{array}{c|cccccccc} \hline
    & t=0 & t=1  & \cdots & t=K-1 & t=K &  \cdots\\
\hline
X(t) & \quad (1,1,\cdots,1) &\quad  (2,2,\cdots, 2)\quad & \cdots & \quad (K,K,\cdots, K) &\quad  (1,1,\cdots,1) \quad   & \cdots\\
\varphi_{n,t} & e_1 & e_2 & \cdots & e_K & e_1 &  \cdots \\ \hline 
\end{array}
\]
The corresponding local distributions are $\varphi_{n,t}=  e_{t+1}\in \R^{1\times K}$ for each $n\in [N]$ and $t\leq K-1$. 
\end{example}

\subsection{Synchronization}
We show first that the system with an irreducible local transition matrix $\bT$ synchronizes if and only if $\bT$ is periodic with period $K$. 
\begin{definition}[Synchronization] We say the system achieves a \emph{synchronization} at time $t_0$ if all sites move identically after $t_0$, i.e., $X_1(t) = X_2(t)=\cdots = X_N(t)$ for all $t \ge t_0$.   
\end{definition}
The following proposition says that if $\bT$ is irreducible and periodic, then the Markov chain $X(\cdot)$ will eventually be periodic. 
\begin{proposition}
    \label{prop:periodicity}
    Suppose that $\bT$ is irreducible and periodic with period $2 \le d \le K$.
    Then
    \begin{enumerate}[label=(\roman*)]
        \item $[K]$ can be decomposed as a finite disjoint union $C_0 \cup \dotsb \cup C_{d-1}$, such that $\bT_{jk} > 0$ only if $j \in C_r$ and $k \in C_{r+1}$ for some $r$ (with the convention $C_d=C_0$). 
        \item For $X(\cdot)$, the collection of states
        $
            S := C_0^{\otimes N} \cup \dotsb \cup C_{d-1}^{\otimes N} \subset [K]^N
        $
        is periodic with period $d$.
        The collection of states $[K]^N \setminus S$ is transient.
    \end{enumerate}
\end{proposition}

\begin{proof}
    (i) This is a standard result; see, e.g.\cite[Theorem 1.8.4]{norris1998markov}. 

    (ii) It follows from part (i) and the global transition matrix in \eqref{eq:T2P} that  
    $$\prob{X(1) \in C_{r+1}^{\otimes N} \mid X(0) \in C_r^{\otimes N}} = 1, \quad \forall\,r\in \{0,1,\dotsc,d-1\},$$ 
    and hence the states in $S$ are periodic with period $kd$ for some (possibly different) $k \in \mathbb{N}$. 
    Now without loss of generality considering $C_0$, for any state $k_0 \in C_0$, by part (i) there exists a chain of states $k_r \in C_r$, $r=1,\dotsc,d-1$, such that $\bT(k_r,k_{r+1})>0$ for each $r=0,\dotsc,d-1$.
    It then follows again from the global transition matrix in \eqref{eq:T2P} that for any state $x_{1:N} \in C_0^{\otimes N}$, one has $\bP^d(x_{1:N},x_{1:N})>0$.
    Therefore, any state in $S$ has period $d$.
    
    It remains to show that the states in $[K]^N \setminus S$ are transient.
    For this, we will show that starting from $X(0)=x_{1:N}$ for any configuration $x_{1:N} \in [K]^N \setminus S$, there is some $t_0 \in \mathbb{N}$ such that $\prob{X(t_0) \in S}>0$. 
    Assume without loss of generality that $x_1=c_0 \in C_0$.
    Since $\bT$ is irreducible, by part (i), there is some $c_1 \in C_1$ such that $\bT_{c_0c_1}>0$.
    Therefore, the nodes in $V_1$ has a positive probability to be $c_1$, i.e., 
    $X(1)$ has a positive probability of taking a configuration in 
    \[
        \{(y_1,\dotsc,y_N) : y_i =c_1, \forall i\in V_1  
        \}, 
    \]
    i.e., the nodes in $V_1= \{1, 2, \ldots,1+n_v, N+1-n_v,\ldots, N\}$ have the same state $c_1$. 
    Again because $\bT$ is irreducible, by part (i), there is some $c_2 \in C_2$ such that $\bT_{c_1c_2}>0$.
    So there is a positive probability that $X(2)$ takes a configuration in
    \[
        \big\{(y_1,\dotsc,y_N) : y_i =c_1, \forall i\in \bigcup_{n \in V_1} V_n   
        \big\},
    \]
    i.e., the nodes in $\bigcup_{n \in V_1} V_n = \{1, 2, \ldots,1+2n_v, N+1-2n_v,\ldots, N\} $ have the same state $c_2$. 
    Continuing in this manner, there is a positive probability for $X$ to jump to some state in $S$ after at most $\lceil \frac{N}{2n_v} \rceil$ steps, i.e., $\prob{X(\lceil \frac{N}{2n_v} \rceil) \in S}>0$. 
\end{proof}
{Note that the proof of Part (ii) only needs the expanding neighborhoods of neighborhoods to cover all the nodes of the graph in finitely many iterations, i.e., for each $i\in [N]$, there exists $t_0<\infty$ such that $\mathcal{N}^{(t_0)}_i = [N]$, where $\mathcal{N}^{(t)}_i := \bigcup_{n \in \mathcal{N}^{(t-1)}_i} V_n$ with $\mathcal{N}^{(0)}_i = V_i$. Thus, Proposition \ref{prop:periodicity} and the rest results of this section can be extended to general connected finite graphs. }

Using Proposition \ref{prop:periodicity}, we have the following characterization of when synchronization occurs.

\begin{theorem}[Synchronization] 
\label{thm:sync}
Suppose that $\bT$ is irreducible. Then, the system will achieve synchronization (with period $K$) if and only if the period of $\bT$ is $K$.
\end{theorem}

\begin{proof} 
    (i) First, we prove the ``if'' direction. Suppose the period of $\bT$ is $K$.
    Then the decomposition of $[K]$ in Proposition \ref{prop:periodicity}(i) must have the form that each $C_r$ is a singleton.
    So the collection of states $S$ in Proposition \ref{prop:periodicity}(ii) is actually just
    \begin{equation}
        \label{eq:sync-states}
        \{(k,\dotsc,k) \in [K]^N : k \in [K]\}.
    \end{equation}
    Also, by Proposition \ref{prop:periodicity}(ii), all states in $[K]^N \setminus S$ are transient and will eventually jump to some state in $S$. 
    Therefore, the system will achieve synchronization with period $K$.

    (ii) Next, we prove the ``only if'' direction. 
    Suppose that the system will achieve synchronization.
    Then there exists some $k_0 \in [K]$ such that, starting from $(k_0,\dotsc,k_0) \in [K]^N$, the system will jump only among the states in \eqref{eq:sync-states}.
    This implies that $\bT_{k_0,k_1}=1$ for some $k_1 \in [K] \setminus \{k_0\}$.
    This further implies, for each $r=1,\dotsc,K-2$, $\bT_{k_r,k_{r+1}}=1$ for some $k_{r+1} \in [K] \setminus \{k_0,\dotsc,k_r\}$, as otherwise we have $k_{r+1}=k_j$ for some $0 \le j \le r$ and it contradicts the assumption that $\bT$ is irreducible.
    Finally, again by irreducibility of $\bT$, we must have $\bT_{k_{K-1},k_0}=1$.
    Therefore, the period of $\bT$ is $K$.
\end{proof}

\begin{remark}
    The arguments in the proof of Proposition {\rm\ref{prop:periodicity}} also reveal that when $\bT$ is reducible, it is still possible that the system will reach a synchronization.
    For example:
    \begin{enumerate}[label=(\roman*)]
    \item 
        If $\bT$ has some transient states, such as $\bT = \begin{bmatrix} 
          0 & 1 & 0 \\
          0 & 0 & 1 \\
          0 & 1 & 0
        \end{bmatrix}$, then eventually, the system will reach a synchronization and oscillate between $(2,\dotsc,2), (3,\dotsc,3) \in [K]^N$.
    \item 
        If $\bT$ has more than one communication class, such as 
        $$\bT = \begin{bmatrix} 
          0 & 1 & 0 & 0 \\
          1 & 0 & 0 & 0 \\
          0 & 0 & 0 & 1 \\
          0 & 0 & 1 & 0
        \end{bmatrix},$$ then eventually the system will reach a synchronization and oscillate only between the configurations $(1,\dotsc,1), (2,\dotsc,2) \in [K]^N$, or only between $(3,\dotsc,3), (4,\dotsc,4) \in [K]^N$.
    \end{enumerate}
\end{remark}

\subsection{Ergodicity}
We show that the system with an irreducible local transition matrix $\bT$ is ergodic if and only if $\bT$ is aperiodic.  

\begin{proposition}\label{prop:irreduciblePT}
    The global transition matrix $\bP$ is irreducible and aperiodic if and only if the stochastic matrix $\bT$ is irreducible and aperiodic.
\end{proposition}

\begin{proof}
    (i) First, we prove the ``if'' direction: 
    Let $\bT$ be irreducible and aperiodic. 
    It suffices to show that there exists $l_0$ such that for all $l\geq l_0$, $\bP^{l}(x,y) >0$ for all $x=(j_1,\dotsc,j_N), y=(k_1,\dotsc,k_N) \in [K]^N$.    
   
    Since $\bT$ is irreducible and aperiodic, there exists some $l_0$ such that for all     
    $l \ge l_0$ and $j,k \in [K]$, $( \bT^{l})_{jk} > 0$. 
    Fix such an $l$. 
    For each $n\in [N]$, $(\bT^{l})_{j_nk_n} > 0$ implies that there exists a sequence $i_1(j_n,k_n),\dotsc,i_{l-1}(j_n,k_n) \in [K]$ such that, writing $k_n=i_l(j_n,k_n)$, $$\bT_{j_ni_1(j_n,k_n)}\bT_{i_1(j_n,k_n)i_2(j_n,k_n)}\dotsm\bT_{i_{l-1}(j_n,k_n)i_l(j_n,k_n)} > 0.$$
    Let $x^{(t)}= (i_t(j_1,k_1), \ldots, i_t(j_N,k_N))$ for $1\leq t\leq l$. 
    Then, noting that $y=x^{(l)}$, we have
    \[
     \bP^l(x,y) \geq  \bP(x,x^{(1)}) \bP(x^{(1)},x^{(2)}) \cdots \bP(x^{(l-1)},x^{(l)}).  
    \]
  
    Meanwhile, since the neighborhood of the node $n$ includes itself, we have $\varphi_n(x_n)\geq \frac{1}{|V_n|}$ for any state $x_{1:N}$. 
    Thus, \eqref{eq:T2P} implies that 
	$\bP(x_{1:N},y_{1:N}) \geq  \prod_{n=1}^N \frac{1}{|V_n|} \bT_{x_n,y_n}$
    for any states $x_{1:N}$ and $y_{1:N}$. Thus, $\bP(x^{(t)},x^{(t+1)}) \geq \prod_{n=1}^N \frac{1}{|V_n|} \bT_{x_n^{(t)},x_n^{(t+1)}} $ for each $t$. 
    Hence, 
    \begin{align*}
        \bP^l(x,y) \ge \prod_{n=1}^N \left( \bT_{j_ni_1(j_n,k_n)}\bT_{i_1(j_n,k_n)i_2(j_n,k_n)}\dotsm\bT_{i_{l-1}(j_n,k_n)i_l(j_n,k_n)}  \cdot \frac{1}{|V_n|^l} \right) > 0.
    \end{align*}
    Therefore $\bP$ is also irreducible and aperiodic.

    (ii) Next, we prove the ``only if'' direction.
    Let $X(\cdot)$ be irreducible and aperiodic.
    
    Suppose $\bT$ is reducible.
    Then there exists $B \subset [K]$ with $B,B^c \ne \emptyset$ such that for all $j \in B, k \in B^c$ and $l \ge 1$, $(\bT^l)_{jk}=0$.
    Then for all $x=(j_1,\dotsc,j_N) \in B^{\otimes N}, y=(k_1,\dotsc,k_N) \in (B^c)^{\otimes N}$ and $l \ge 1$, $\bP^{l}(x,y)=0$. 
    Therefore $\bP$ is reducible.
    This leads to a contradiction, and hence $\bT$ must be irreducible.

    Suppose $\bT$ is irreducible but periodic with period $d\geq 2$. 
    Then, by Proposition \ref{prop:periodicity}(ii), $X(\cdot)$ is not irreducible, which is a contradiction.
    This completes the proof. 
\end{proof}

We conclude in the next theorem that the Markov chain is exponentially ergodic when the local transition matrix is irreducible and aperiodic. Its proof is postponed to Section \ref{sec:prelimMC}.  
\begin{theorem}
    \label{thm:ergodic}
    Assume that the stochastic matrix $\bT$ is irreducible and aperiodic. Then, $X(\cdot)$ is irreducible, aperiodic, and hence ergodic.
    In particular, there exists a unique stationary distribution $\pi$, all states are positive recurrent, $\lim_{t \to \infty} \prob{X(t) = \cdot} = \pi(\cdot)$, and there exist $C \in (0,\infty)$, $\rho \in (0,1)$ and $t_0 \in \mathbb{N}$ such that
    \begin{equation}
        \label{eq:geometic-cvg}
        \|\bP^t(x,\cdot)-\pi\|_{TV} \le C\rho^t, \quad \forall\,t \ge t_0, \: x \in [K]^N,
    \end{equation}
    where $\|\cdot\|_{TV} := \|\cdot\|_1$ denotes the total variation norm. 
\end{theorem}

The following corollary is a particular case of Theorem \ref{thm:ergodic}, with $t_0=1$ in its proof.

\begin{corollary}
    Suppose $\bT_{jk}>0$ for all $j,k \in [K]$.
    Then all the statements in Theorem {\rm\ref{thm:ergodic}} hold.
    In fact, there exists $\rho \in (0,1)$ such that
    \[
        \|\bP^t(x,\cdot)-\pi\|_{TV} \le 2\rho^t, \quad \forall\,t \ge 1, \: x \in [K]^N.
    \]
\end{corollary}

\begin{remark}
    Proposition {\rm\ref{prop:periodicity}}, Theorem {\rm\ref{thm:sync}}, and Theorem {\rm\ref{thm:ergodic}} now fully characterize the long-time behavior of the system when the local transition matrix $\bT$ is irreducible.
    In particular, Proposition {\rm\ref{prop:periodicity}} shows that if $\bT$ is periodic with period $d < K$, then the system will neither achieve a synchronization (as at least one class of $C_0,\dotsc,C_{d-1}$ contains more than one state) nor is exponentially ergodic.
    It then follows from Theorems {\rm\ref{thm:sync}} and {\rm\ref{thm:ergodic}} that the system will achieve a synchronization if and only if $\bT$ is periodic with period $K$, while the system is exponentially ergodic if and only if $\bT$ is aperiodic.
\end{remark}


\section{Global and local transition matrices}\label{sec:PT}
We study in this section the relation between the global transition matrix $\bP \in \R^{K^N\times K^N}$ and the local transition matrix $\bT\in \R^{K\times K}$, and their associated invariant measures. 
Throughout this section, we assume that $\bT$ is irreducible and aperiodic, and so is $\bP$ by Proposition \ref{prop:irreduciblePT}. 

\subsection{Local and global transition matrices and associated invariant measures} \label{sec:invTP}
We show that there exists a 1-1 map between the global and local transition matrices $\bP$ and $\bT$. Furthermore, we show that $\pitilde$ is the marginal of $\pi$, where $\pitilde \in \R^{K}$ and $\pi\in \R^{K^N}$ denote the unique invariant measures of $\bP$ and $\bT$, respectively. That is, the marginal distribution of the Markov chain's invariant measure is the same as the invariant measure of $\bT$.  
However, the marginal distribution does not necessarily determine $\bT$, as we show in Example \ref{ex:twoT_diff_inv} that there are two $\bT$'s leading to the same $\pitilde$ and different $\pi$'s; nor does the joint invariant distribution, as we show in Example \ref{ex:1pi-2T} that there are two $\bT$'s leading to the same $\pitilde$ and $\pi$.  

\begin{proposition}[Shift-invariance]
    \label{prop:shift-invariance}
    The transition matrix $\bP$ is shift invariant. Therefore, the invariant measure $\pi$ is shift-invariant, namely, for each $(x_1,\dotsc,x_N) \in [K]^N$,
    $$\pi(x_1,\dotsc,x_N) = \pi(x_2,\dotsc,x_N,x_1).$$
\end{proposition}

\begin{proof}
    Recall the neighborhood $V_n= \{n-n_v,\ldots, n, n+1,\ldots, n+n_v\}$.
    Clearly, the graph $G=(V,E)$ is invariant under the shift of the node indices $\{1,\dotsc,n\}$.
    As the transition of $X(\cdot)$ in \eqref{eq:MC} and \eqref{eq:PT} depends on states through the empirical distribution of the neighborhood, the transition matrix $\bP$ is shift invariant.
    Since the invariant measure $\pi$ is unique, we also have shift-invariance of $\pi$.
\end{proof}

The shift-invariance originates from the symmetry of the N-cycle graph and the $n_v$-distance neighborhoods $V_n$. One can generalize the above proposition and the next theorem to general symmetric or vertex-transitive finite graphs.   

The next theorem shows that there is a 1-1 correspondence between the local and global transition matrices and that $\pitilde$ is the marginal of $\pi$. 

\begin{theorem}[1-1 map between $\bP$ and $\bT$, $\tilde{\pi}$ as the marginal of $\pi$]
\label{thm:1-1PT}
There is a 1-1 map between $\bT\in \R^{K\times K}$ and $\bP\in \R^{K^N\times K^N}$.
Furthermore,  
denote by $\pi_n$ the $n$-th marginal distribution of $\pi$, i.e., $\pi_n(k) = \lim_{t \to \infty} \mathbb{P}(X_n(t)=k)$, $\forall\,k =1,\dotsc,K$. Then $\pi_n = \pitilde$.
\end{theorem}

\begin{proof} 
Eq.\eqref{eq:T2P} shows that $\bT$ uniquely determines $\bP$. It also implies that 
\begin{equation*}
\bT_{jk} = \bP(x_{1:N},y_{1:N})^{1/N}, \, \text{ with } x_{1:N} =(j,\dotsc,j), \, y_{1:N}=(k,\dotsc,k).
\end{equation*}
Therefore there is a 1-1 map between $\bT\in \R^{K\times K}$ and $\bP\in \R^{K^N\times K^N}$.

Next, we prove that the marginal distribution $\pi_n$ of the invariant measure $\pi$ of $\bP$ is the same as the invariant measure $\tilde{\pi}$ of $\bT$. 
    By shift-invariance in Proposition \ref{prop:shift-invariance}, $\pi_n = \pi_m$ for all $n,m\in [N]$.
    Since $\pi$ is the invariant measure of $\bP$,  we have $\pi= \pi \bP$. Applying \eqref{eq:T2P}, we can write 
    $$
    \pi(y_{1:N}) = \sum_{x_{1:N} \in [K]^N} \pi(x_{1:N}) \prod_{n\in [N]} \frac{1}{|V_n|} \sum_{n'\in V_n} \bT_{x_{n'}y_n}.$$
    Note that $\frac{1}{|V_n|} \sum_{n'\in V_n} \bT_{x_{n'}y_n}$ depends on $y_{1:N}$ only through $y_n$.
    Summing over $y_{2:N} \in [K]^{N-1}$, we have
    $$
    \pi_1(y_1) = \sum_{x_{1:N} \in [K]^N} \pi(x_{1:N}) \frac{1}{|V_1|} \sum_{n'\in V_1} \bT_{x_{n'}y_1} = \frac{1}{|V_1|} \sum_{n'\in V_1} \sum_{x_{n'} \in [K]} \pi_{n'}(x_{n'}) \bT_{x_{n'}y_1}.$$
    By symmetry of $\pi$, we have $\pi_1=\pitilde$ since
    $\pi_1(y_1)=\sum_{k \in [K]} \pi_1(k)\bT_{ky_1}\,.$
\end{proof}

However, $\pi$ is not a product measure of $\pitilde$ and it is not completely determined by $\pitilde$. 
The following example shows that multiple systems with different $\bT$'s can have different $\pi$'s but the same $\pitilde$.
\begin{example}[Same marginal, different invariant measures, $\bT$ and $\bP$'s] 
\label{ex:twoT_diff_inv}
    Consider $K=2$, $N=2$ and $V_n=\{1,2\}$ without doubling counting the other vertex. Each of the following $\bT$ has $\pitilde = \begin{bmatrix}  2/3,  1/3 \end{bmatrix}$ as invariant measure. But they have different $\bP$'s and $\pi$'s.  
    \begin{align*}
    \text{ (i). }\quad 	\bT&=\begin{bmatrix}
                              1/2 & 1/2 \\
                               1 & 0
                              \end{bmatrix}, \quad 
          \pi      = \begin{bmatrix}
                             10/21, 4/21, 4/21, 3/21
                     \end{bmatrix}, \quad  
          \bP   =\begin{bmatrix}
              1/4 & 1/4 & 1/4 & 1/4 \\
              9/16 & 3/16 & 3/16 & 1/16 \\
              9/16 & 3/16 & 3/16 & 1/16 \\
              1 & 0 & 0 & 0
              \end{bmatrix}; \\
      \text{ (ii).} \quad \bT &=\begin{bmatrix}
                             3/4 & 1/4 \\  1/2 & 1/2 
                            \end{bmatrix}, \quad  
         \pi      = \begin{bmatrix}
                             14/31 , 20/93, 20/93 , 11/93
                     \end{bmatrix}, 
          \bP=\begin{bmatrix}
            9/16 & 3/16 & 3/16 & 1/16 \\
            25/64 & 15/64 & 15/64 & 9/64 \\
            25/64 & 15/64 & 15/64 & 9/64 \\
            1/4 & 1/4 & 1/4 & 1/4
           \end{bmatrix};\\
       \text{ (iii).} \quad \bT &=\begin{bmatrix}
                             2/3 & 1/3 \\  2/3 & 1/3 
                            \end{bmatrix}, \quad  
         \pi      = \begin{bmatrix}
                             4/9 , 2/9, 2/9 , 1/9
                     \end{bmatrix},\quad 
        \bP=\begin{bmatrix}
            4/9 &  2/9 &  2/9 &  1/9 \\
            4/9 &  2/9 &  2/9 &  1/9 \\
            4/9 &  2/9 &  2/9 &  1/9 \\
            4/9 &  2/9 &  2/9 &  1/9
        \end{bmatrix}. 
    \end{align*}
\end{example}

An interesting question is whether a $1$--$1$ map exists between $\pi$ and $\bP$ for the PCA. 
Clearly, $\bP$ uniquely determines $\pi$. The other direction is not true for a general Markov chain: there can be multiple transition matrices with the same invariant measure (see e.g.\ Example \ref{ex:twoT_diff_inv} where two different $\bT$'s lead to the same invariant measure $\widetilde{\pi}$). 
However, given the special structure of $\bP$ being determined by a local transition matrix $\bT\in \R^{K\times K}$, which has only $K(K-1)$ ``free entries'', one may question whether the invariant measure $\pi\in \R^{K^N}$ can uniquely determine $\bT$ and consequently $\bP$. 
The following example demonstrates that the answer is no. 

\begin{example}[Same invariant measure, different $\bT$ and $\bP$'s]
\label{ex:1pi-2T} We still have 
  the same $\pi$ (and hence $\tilde{\pi}$) as in Example {\rm \ref{ex:twoT_diff_inv}}(ii) if $\bT=\begin{bmatrix}
        7/12 & 5/12 \\
        5/6 & 1/6
    \end{bmatrix}$, which leads to 
\[
 \bP=\begin{bmatrix}
        49/144 & 35/144 & 35/144 & 25/144 \\
        289/576 & 119/576 & 119/576 & 49/576 \\
        289/576 & 119/576 & 119/576 & 49/576 \\
        25/36 & 5/36& 5/36 & 1/36
    \end{bmatrix}. 
\] 
\end{example}

\subsection{Lipschitz dependence on the local transition matrix}
We show the Lipschitz dependence of the global transition matrix and its invariant measure on the local transition matrix.
For a stochastic matrix $\bA$, let $\tau(\bA) := \frac{1}{2} \max_{x,x'} \|\bA(x,\cdot) - \bA(x',\cdot)\|_1$.

\begin{theorem}[Lipschitz dependence on $\bT$]
\label{thm:Lip_P_T}
Given $\bT^{(1)}$ and $\bT^{(2)}$, let $\bP_1, \bP_2$ be the corresponding global transition matrices and $\pi_1,\pi_2$ be the stationary measures. 
Then  
\begin{align}
        \| \bP_1 - \bP_2 \|_1 &\le N K^{N-1} \min\{\|\bT^{(1)}-\bT^{(2)}\|_1, K \|\bT^{(1)}-\bT^{(2)}\|_{1,1}\},         \label{eq:DeltaP-L1} \\
        \|\bP_1 - \bP_2 \|_2 &\le N(KC_K)^{N/2}  \|\bT^{(1)}-\bT^{(2)}\|_2,         \label{eq:DeltaP-L2}
 \end{align}
where $\|\bA\|_p= (\sum_{ij} |\bA_{i,j}|^p)^{1/p}$ is the $\ell_p$ norm and $\|\bA\|_{1,1} = \sup_{\|\varphi\|_1 \le 1} \|\varphi\bA\|_1$ is the $1$-operator norm for a matrix $\bA$,
   and  
    $
    C_K := \max\left\{\max_{j \in [K]} \sum_{k \in [K]} (\bT^{(1)}_{j,k})^2, \max_{j \in [K]} \sum_{k \in [K]} (\bT^{(2)}_{j,k})^2\right\} \leq 1.
    $

Also, there exists $\ell_0 \in \mathbb{N}$ such that $\tau(\bP_1^{\ell_0})<1$ and
\begin{align*}
\|\pi_1- \pi_2 \|_1&\leq \frac{1+(\ell_0-1)K^N}{1-\tau(\bP_1^{\ell_0})} \|\bP_1-\bP_2\|_{1} \\
& \leq \frac{1+(\ell_0-1)K^N}{1-\tau(\bP_1^{\ell_0})} N K^{N-1}\min\{\|\bT^{(1)}-\bT^{(2)}\|_1, K \|\bT^{(1)}-\bT^{(2)}\|_{1,1}\}.	
\end{align*}
\end{theorem}

\begin{figure}
	\centering
	\includegraphics[width=0.45\textwidth]{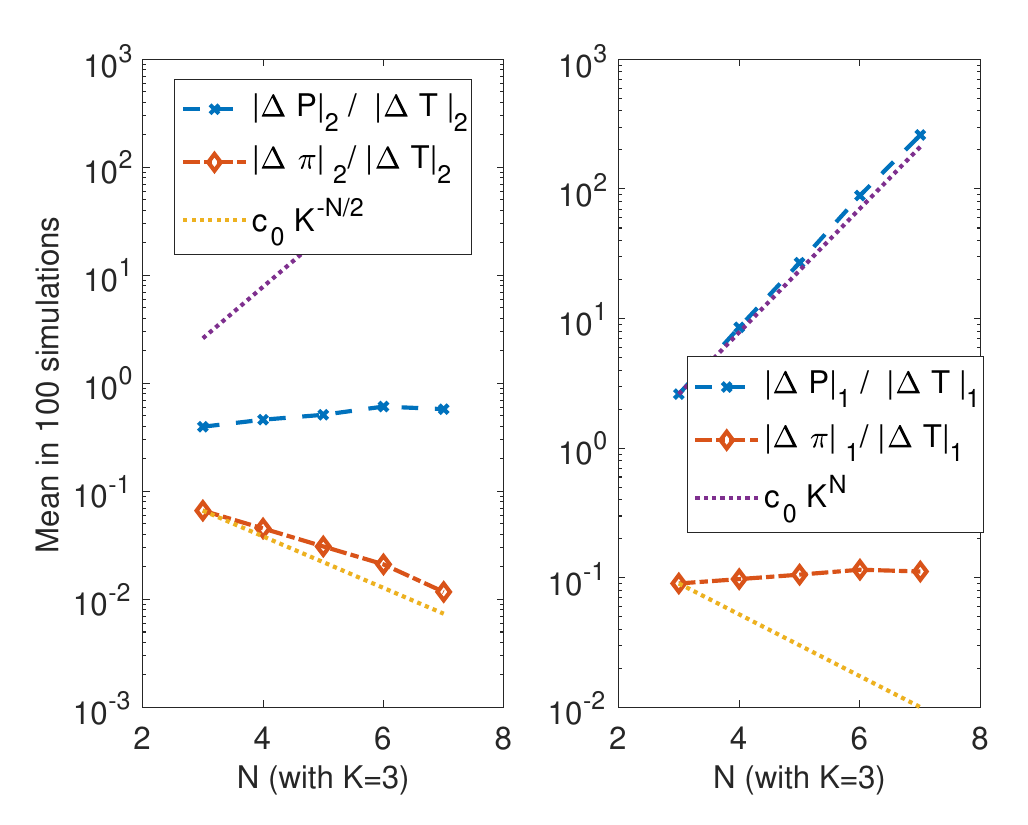}
	\includegraphics[width=0.45\textwidth]{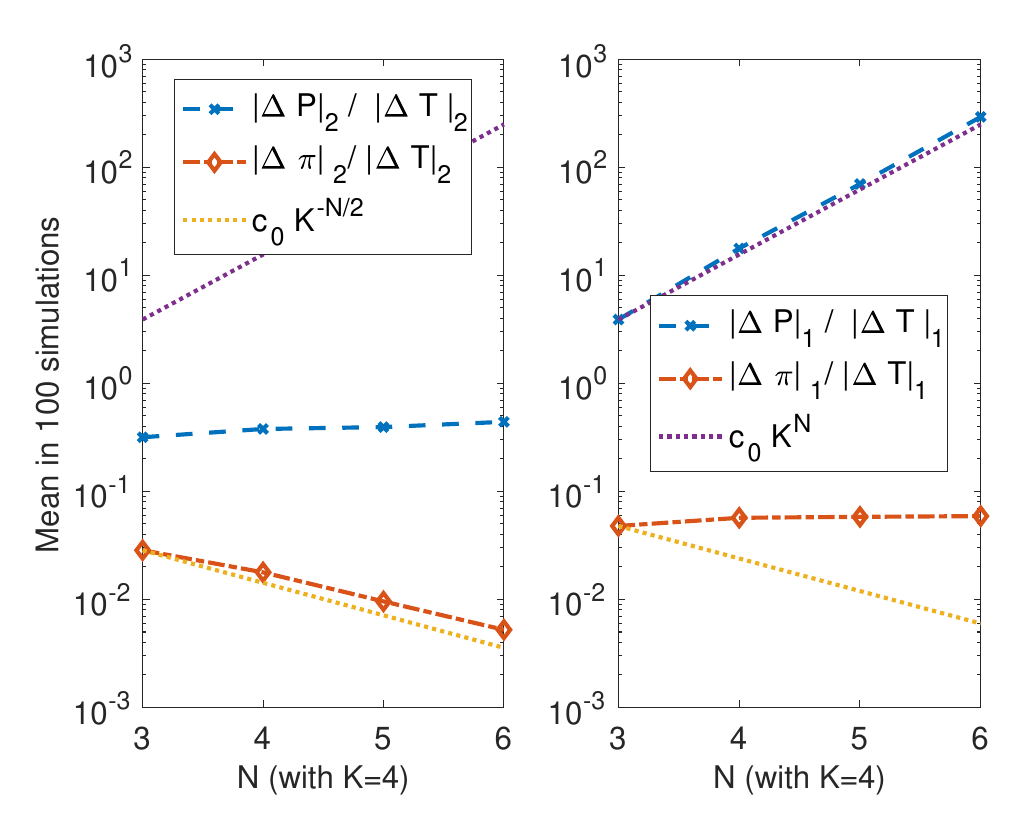}
	\caption{ Mean ratios  $\frac{\|\Delta \bP\|_{p}}{\|\Delta \bT\|_{p} }$ and $\frac{\|\Delta \pi\|_{p}}{\|\Delta \bT\|_{p} }$ with $p=1,2$, where each $(N,K)$-pair is computed using 100 random $\bT$ with entries sampled from uniform [0,1] followed by a row-normalization. Here the neighborhood has $n_v= \min\{3,\floor{N/2}\}$.  
	Note that $\frac{\|\Delta \bP\|_{1}}{\|\Delta \bT\|_{1} } = O(K^N)$, agreeing with Theorem \ref{thm:Lip_P_T}. Also, note that $\frac{\|\Delta \pi\|_{1}}{\|\Delta \bT\|_{1} } = O(1)$, $\frac{\|\Delta \bP\|_{2}}{\|\Delta \bT\|_{2} } = O(1)$, and $\frac{\|\Delta \pi\|_{2}}{\|\Delta \bT\|_{2} } = O(K^{-N/2})$.  
	  }
	\label{Fig:Lip_PT_ratio}
\end{figure}

We postpone its proof to Section \ref{sec:Proof-Lip_P_T} in the appendix but comment on the optimality of the upper bounds here. 
The above Lipschitz dependence on $\bT$ for the invariant measure is local because the constant $\frac{1}{1-\tau(\bP_1)}$ can depend on $\bT^{(1)}$. 
    Numerical tests, as shown in Figure {\rm \ref{Fig:Lip_PT_ratio}}, suggest that the bound for $\bP$ is optimal; however the bound for $\pi$ may not be optimal. 
Also, the bounds in Frobenius norm in \eqref{eq:DeltaP-L2} is suboptimal since the constant $N(KC_K)^{N/2}$ may increase exponentially in $N$.  
   Intuitively, if $\tau(\bP_1)<1$ and the entries of $\bP_1 - \bP_2$ are of a similar order, then using \eqref{eq:DeltaP-L1} of  Theorem \ref{thm:Lip_P_T} one would expect
    \[
    \| \bP_1 - \bP_2 \|_2 \approx C K^{-N}\| \bP_1 - \bP_2 \|_1 \le CNK^{-1} \|\bT^{(1)}-\bT^{(2)}\|_1 \le CN \|\bT^{(1)}-\bT^{(2)}\|_2, 
    \]
    which agrees with Figure {\rm\ref{Fig:Lip_PT_ratio}}.  

\section{Inference of the local transition matrix}\label{sec:inference}
Inference is the first step in the application of the PCA model. In this section, we construct least squares estimators (LSEs) to infer the local transition matrix from various types of data. Data may consist of multiple trajectories, a long trajectory, or ensemble sequences without trajectory information, discussed in Sections \ref{sec:lse_traj}--\ref{sec:lse_ens}. For each of these cases, we specify identifiability conditions and prove that the estimators are asymptotically normal. Additionally, we show $\bT$ is non-identifiable from the stationary distribution in Section \ref{sec:infer_from_inv}. Furthermore, we establish non-asymptotic bounds for the estimators in Section \ref{sec:non-asymp}.  Finally, the numerical tests in Section \ref{sec:num} demonstrate that the LSE with trajectory information is more accurate than the LSE without trajectory information, while both converge at the rate $M^{-1/2}$, in agreement with theory. 

\subsection{LSE from multiple trajectories }\label{sec:lse_traj} 
Consider first the inference of $\bT\in \R^{K\times K}$ from data consisting of $M$ independent trajectories:
\[
\textbf{Trajectory Data:} \quad \quad \{X^m(t), t=1,\ldots, L\}_{m=1}^M.
\]
We estimate each column of $\bT$ through least squares, followed by a normalization. 
The least squares estimator minimizes the loss function 
\begin{equation}\label{eq:lossFn-traj}
	\widehat {\bT}_M = \argmin{\bT \in \R^{K\times K}}  \mathcal{E}(\bT) ,\quad \text{ with }  \mathcal{E}(\bT) : = \frac1M\frac1L\frac1N\sum_{m=1}^M\sum_{t=1}^L \sum_{n=1}^N \sum_{k=1}^K \big | c_{n,t}^m(k)  -  \varphi_{n,t-1}^m \bT_{\cdot, k}\big|^2\,,  
\end{equation} 
where $c_{n,t}^m \in \R^{1\times K}$ denotes the empirical distribution $X_n^m(t)$:  
\begin{equation}\label{eq:c_nl}
c_{n,t}^m = ( c_{n,t}^{m}(1),\ldots,  c_{n,t}^{m}(K) ) \in \R^{1\times K} \quad \text{ with }   c_{n,t}^{m}(k) := \delta_{X_n^m(t)}(k), 
\end{equation}
and $\varphi_{n,t}^m\in \R^{1\times K}$ is the empirical distribution of the states of sample $X^m(t)$ in the neighborhood $V_n$ of the vertex $n$, as defined in \eqref{eq:phi}. By solving the zero of the gradient of this loss function with respect to each $\bT_{\cdot, k}$, we obtain the least squares estimator from a system of normal equations $K$ with a shared normal matrix $\bA_{M}\in \R^{K\times K}$: 
 \begin{equation}\label{eq:traj-LSE} 
\begin{aligned}
\widehat {\bT}_M( \cdot,k) &= \bA_{M} ^{\dag}  \bb_{M}(\cdot, k), \quad  1\leq k\leq K,  
  \\
\bA_{M} & =\frac1{MLN}\sum_{t,n,m=1}^{L,N,M}  ( \varphi_{n,t-1}^{m})^\top  \varphi_{n,t-1}^{m}\,, \quad 
\bb_{M}(\cdot,k)  =  \frac1{MLN}\sum_{t,n,m=1}^{L,N,M} \left(\varphi_{n,t-1}^{m}\right)^\top  c_{n,t}^{m}(k). 
\end{aligned}
\end{equation}
Here, $ \bA_{M} ^{\dag}$ denotes the Moore-Penrose pseudo inverse of $\bA_M$. 

In practice, instead of using pseudo-inverse, we obtain $\widehat \bT_M$ by using least squares with non-negative constraints and then row-normalize the solution. The non-negative constraints help to avoid possible negative entries caused by sampling error.

\paragraph{Identifiability from the large sample limit.}
To analyze the estimator, we first examine the inference problem in the large sample limit. Denote the large sample limit normal matrix and normal vectors by 
\begin{equation}\label{eq:Abinfty}
\begin{aligned}	
	\bA_\infty  =  \frac1{LN} \sum_{t=1}^L \sum_{n=1}^N   \E[ \varphi_{n,t-1}^\top  \varphi_{n,t-1} ], \quad  
	\bb_\infty(\cdot, k)  =  \frac1{LN} \sum_{t=1}^L \sum_{n=1}^N   \E[ \varphi_{n,t-1}^\top c_{n,t}(k)], \,  1\leq k\leq K. 
\end{aligned}
\end{equation}
\begin{assumption}[Identifiability condition: multi-trajectory data]
\label{assump:Ainvertible} 
The distribution of the data satisfies that the matrix $\bA_\infty $ in \eqref{eq:Abinfty} is non-singular. 
\end{assumption}

The assumption \ref{assump:Ainvertible} holds in general, as long as the random vectors $\{\varphi_{n,t}\}$ span $\R^K$ with positive probability, as the next lemma shows.  
\begin{lemma} 
Assumption {\rm\ref{assump:Ainvertible}} holds except when there exists a fixed vector $v\in\R^{1\times K}$ such that $v \varphi_{n,t-1}^\top=0$ a.s.~for all $n$ and $l$. In particular, when $\bT$ is irreducible and aperiodic, Assumption {\rm\ref{assump:Ainvertible}} holds as long as $L\geq l_0: = \argmin{t\geq 0} \{t: \bP^t(x,y)>0, \forall x,y\in [K]^N\}$.  
\end{lemma}
\begin{proof} Recall that the covariance matrix $\E[Z^\top Z]$ of a random vector $Z\in \R^{1\times K}$ is singular iff there exists a vector $v\in \R^{1\times K}$ such that $ Z v^\top=0$ a.s., which is true because $0=v \E[Z^\top Z]v^\top = \E[| Zv^\top |^2]$ iff $ Zv^\top=0$ a.s.. Applying this fact to the random vectors $ \varphi_{n,t-1} $, we find that $\bA_\infty$ is singular iff there exists a fixed vector $v$ such that $\varphi_{n,t-1} v^\top=0$ a.s.~for all $n$ and $t$. 

When $\bT$ is irreducible and aperiodic, $\bP$ is also irreducible and aperiodic by Proposition \ref{prop:irreduciblePT}.  The number $l_0$ is well-defined and is finite. Then, regardless of what the initial condition is, the states $\{X(t)\}_{t\leq l_0}$ visit all possible states with a positive probability, so $\{\varphi_{n,t}\}_{t\leq l_0}$ span $\R^K$ with a positive probability. Hence, one cannot find a $v$ such that $\varphi_{n,t_0} v^\top=0$ a.s.~for all $n$. 
\end{proof}
 
The exceptions are extreme. For the system in either Example \ref{exp:T_move2next} or Example \ref{exp:T-permutation}, the normal matrix $\bA_\infty$ is singular when $L\leq K-2$, and is non-singular once $L\geq K-1$. For the system in Example \ref{exp:T_move2next}, since  $\varphi_{n,t}= e_{t+1}\in \R^{K\times 1}
 $  for each $n\in [N]$ and $0\leq t\leq K-2$, we can take $v=(0,\dotsc,0,1)= e_K$ so that  $v^\top\varphi_{n,t}=0$ for all $n$ and $t \le K-2$. In this case, $\bA_\infty$ is singular for any $L \le K-2$. On the other hand, if $L\geq K$, the resulting normal matrix $\bA_\infty$ becomes full rank.

\begin{lemma}\label{lemma:Ainv} 
Under Assumption {\rm\ref{assump:Ainvertible}}, $\bT = \bA_\infty^{-1}\bb_\infty$ with $\bA_\infty$ and $\bb_\infty$ in \eqref{eq:Abinfty}.  
\end{lemma}
\begin{proof}
For each $1\leq k\leq K$,  recall that $\varphi_{n,t-1}(k) = \frac{1}{|V_n|} \sum_{j\in V_n} \delta_{X_j(t-1)}(k)$ in \eqref{eq:phi} and $ c_{n,t}(k)= \delta_{X_n (t)}(k)$ in \eqref{eq:c_nl}. 
Then 
\begin{align*}
\E[ \varphi_{n,t-1}^\top c_{n,t}(k) ]  & = \E\left[  \varphi_{n,t-1}^\top \E[ \delta_{X_n(t)}(k) \mid X(t-1)] \right]  
 =  \E\left[   \varphi_{n,t-1}^\top  \varphi_{n,t-1}  \bT_{\cdot,k}\right]
\end{align*}
using the fact that
$\E[ \delta_{X_n(t)}(k) \mid X(t-1)] = \prob{ X_n(t) = k | X(t-1) } 
 =    \varphi_{n,t-1} \bT_{\cdot,k} $ by \eqref{eq:PT}. Hence, 
\begin{align*}
\bb_\infty(\cdot, k)  =  \frac1{LN} \sum_{t=1}^L \sum_{n=1}^N   \E[\varphi_{n,t-1}^\top  c_{n,t}(k)]= \frac1{LN} \sum_{t=1}^L \sum_{n=1}^N  \E\left[ \varphi_{n,t-1}^\top \varphi_{n,t-1}\right] \bT_{\cdot,k}  . 
\end{align*}
In other words, 
$\bb_\infty =  \frac1{LN} \sum_{t=1}^L \sum_{n=1}^N  \E\left[\varphi_{n,t-1}^\top  \varphi_{n,t-1}  \right] \bT = \bA_\infty  \bT.
$
Therefore, $\bT = \bA_\infty^{-1}\bb_\infty$, where $\bA_\infty$ is invertible by Assumption {\rm\ref{assump:Ainvertible}}. 
\end{proof}

\paragraph{Asymptotic normality.} We show next that the LSE is asymptotically normal. 

\begin{theorem} \label{thm:AN_lse} 
Under Assumption {\rm\ref{assump:Ainvertible}}, for each $ k\in [K]$, the estimator $\widehat \bT_M(\cdot,k)$ in \eqref{eq:traj-LSE} is asymptotically normal;  that is, 
 \begin{equation}\label{eq:LSE_AN}
\sqrt{M} \left( \widehat \bT_M(\cdot,k) -  \bT(\cdot,k) \right) \to  \mathcal{N}(0, \bA_\infty^{-1} \Sigma_k \bA_\infty^{-1} ), 
 \end{equation}
 where $\Sigma_k\in \R^{K\times K}$ is the average covariance  of $ \widetilde{c}_{n,t} := c_{n,t}(k)\varphi_{n,t-1}- \E[ c_{n,t}(k)\varphi_{n,t-1}]$: 
\begin{equation}\label{eq:Sigma_k}
\begin{aligned}
\Sigma_k 
=  \frac1{L^2N^2} \sum_{t, t'=1}^L \sum_{n,n'=1}^N \E\left[ \widetilde{c}_{n,t}  \widetilde{c}_{n,t}^\top\right]. 
\end{aligned}
\end{equation}
\end{theorem}
\begin{proof} 
For $1\leq k\leq K$, denote 
\begin{equation}\label{eq:Ab-m}
\begin{aligned}
\bA_{L,N}^m: =\frac1{LN}\sum_{t,n}^{L,N}  ( \varphi_{n,t-1}^{m})^\top   \varphi_{n,t-1}^{m}\,,   \quad 
\bb_{L,N}^m(\cdot,k): &=\frac1{LN}\sum_{t,n}^{L,N} c_{n,t}^m(k)   \varphi_{n,t-1}^{m}\,, 
\end{aligned}	
\end{equation}
where $c_{n,t}(k)$ is defined in \eqref{eq:c_nl}. 
 Note that $\{\bA_{L,N}^m\}_{m=1}^M$ and $\{\bb_{L,N}^m(\cdot,k)\}_{m=1}^M$ are independently identically distributed, and $\bA_{M} = \frac{1}{M}\sum_{m=1}^M \bA_{L,N}^m$ and $\bb_{M}(\cdot,k)= \frac{1}{M}\sum_{m=1}^M \bb_{L,N}^m(\cdot,k)$. Thus, by the strong law of large numbers, 
$\bA_{M}\xrightarrow{a.s.} \bA_{\infty}
$; and by the central limit theorem, 
$$\sqrt{M} ( \bb_{M}(\cdot,k)-\bb_\infty(\cdot,k)) \xrightarrow{d} \mathcal{N}(0,\Sigma_k )$$ 
for each $k$, where the matrix $\Sigma_k$ in \eqref{eq:Sigma_k} is the covariance matrix of $\bb_{L,N}^m(\cdot,k)$:  
\begin{align*}
\Sigma_k &= \E[  \left( b_{L,N}^m(\cdot,k)-\E[b_{L,N}^m(\cdot,k)] \right)   \left( b_{L,N}^m(\cdot,k)-\E[b_{L,N}^m(\cdot,k)] \right)^\top  ].  
\end{align*}

Then, Lemma \ref{lemma:Ainv_b_AN} implies $\sqrt{M} \big(\bA_M^\dag \bb_M(\cdot,k) - \bA_\infty^{-1} \bb_\infty(\cdot,k) \big) \xrightarrow{d} \mathcal{N}(0, \bA_\infty^{-1} \Sigma_k \bA_\infty^{-1} )$. The asymptotic normality of the LSE in \eqref{eq:LSE_AN} follows from $\bT_M(\cdot,k)= \bA_M^\dag \bb_M(\cdot,k)$ in \eqref{eq:traj-LSE} and Lemma \ref{lemma:Ainv}. 
\end{proof}

The following lemma is a slight extension of Slutsky's theorem; we will use it repeatedly to study the asymptotic normality of least squares estimators. Its proof is included for completeness.  
\begin{lemma}\label{lemma:Ainv_b_AN}
Suppose that  $A_M \xrightarrow{a.s.} A$ and $\sqrt{M} \big(b_M - b\big) \xrightarrow{d} \mathcal{N}(0, B)$ as $M \to \infty$, where $A, B$ are two symmetric strictly positive definite matrices.  Then $A_M^{\dag}b_M$ is asymptotically normal, i.e., $\sqrt{M} \big( A_M^{\dag}b_M - A^{-1} b\big) \xrightarrow{d} \mathcal{N}(0,A^{-1}BA^{-1}) $.
\end{lemma}
\begin{proof}
First, we have $A_M^\dag \xrightarrow{a.s.} A^{-1}$ since $A$ is invertible and  $A_M \xrightarrow{a.s.} A$. Specifically, the almost sure convergence of $A_M$ implies that $\prob{\lim_{M\to\infty} \|A_M- A \|_F =0 } = 1 $. Then, Weyl's inequality $|\lambda_{min}(A_M)-  \lambda_{min}(A)|  \leq \|A_M-A\|_{op}$, where $\lambda_{min}(A)$ denotes the minimal eigenvalue of $A$,  implies that  $\prob{ \lim_{M\to \infty}\lambda_{min}(A_M) = \lambda_{min}(A) >0} =1$. Thus, $A_M^\dag \xrightarrow{a.s.}  A^{-1}$.

Next, combining $A_M^\dag \xrightarrow{a.s.} A^{-1}$ with $\sqrt{M}( b_M - b)\xrightarrow{d} \mathcal{N}(0,B)$,  we have, by Slutsky's theorem, $\sqrt{M} \big( A_M^{\dag}(b_M - b) \big) \xrightarrow{d} A^{-1}\mathcal{N}(0, B)= \mathcal{N}(0, A^{-1}B A^{-1})$ and $A_M^\dag b \xrightarrow{a.s.} A^{-1}b$. 
Therefore, using $A_M^{\dag}b_M = A_M^{\dag}(b_M - b) + A_M^{\dag}b$, we obtain $\sqrt{M} \big(A_M^{\dag}b_M - A^{-1}b\big) \xrightarrow{d} \mathcal{N}(0, A^{-1}BA^{-1})$.
\end{proof}


\paragraph{Comparison with the maximal likelihood estimator (MLE).} 
The maximal likelihood estimator involves a constrained non-convex optimization, making it computationally more costly and theoretically more complex to analyze than the above LSE. 
The MLE maximizes the likelihood of data:  
\begin{equation}\label{e:log-MLE2}
\begin{aligned}
\widehat \bT &= \argmin{\bT \text{ stochastic} 
} 
\mathcal{E}_{mle}(\bT),  \quad  \text{ where }\, 
 \mathcal{E}_{mle}(\bT) := \frac{1}{M}\frac{1}{L}\frac{1}{N} \sum_{m=1}^M\sum_{t=1}^L \sum_{n=1}^N - \log (  \varphi_{n,t-1}^m \bT_{\cdot,X_n^m(t)} ) . 
\end{aligned}
\end{equation}
 Note that, unlike the LSE, it is necessary to consider the optimization with respect to stochastic matrices $\bT$, because otherwise, the likelihood has no maximum. The uniqueness of the minimizer for the constrained optimization of such a nonconvex function is relatively complicated to analyze. The asymptotic normality of the MLE is nontrivial since it relies on uniqueness. Also, while optimization algorithms can easily compute a minimizer, it remains difficult to provide a performance guarantee.

\subsection{LSE from a single long trajectory for ergodic systems}

Suppose that the system is ergodic, and we estimate the local transition matrix from
\[
\textbf{A long trajectory Data:} \quad \quad \{X(t), t=0,\ldots, L\}.
\]
Under the ergodicity assumption, the estimation is the same as the previous case with $M=1$, and we define the estimator by 
 \begin{equation}\label{eq:Ab_k_1traj}
\begin{aligned}
\widehat {\bT}_L( \cdot,k) &= \bA_{L} ^{\dag}  \bb_{L}(\cdot, k), \quad 1\leq k\leq K,   
  \\
\bA_{L} & =\frac1{LN}\sum_{t,n=1}^{L,N}   \varphi_{n,t-1}^\top  \varphi_{n,t-1}\,, \quad 
\bb_{L}(\cdot,k)  =  \frac1{LN}\sum_{t,n=1}^{L,N} \left(\varphi_{n,t-1}\right)^\top  c_{n,t}(k).   
\end{aligned}
\end{equation}
Here $ \bA_{L} ^{\dag}$ denotes the Moore-Penrose pseudo-inverse of $\bA_L$.

Similarly to the previous section, the large sample limit helps us specify the identifiability condition. Denote the large sample limit normal matrix and normal vectors by 
\begin{equation}\label{eq:Abinfty_1traj}
\begin{aligned}	
	\bA_\infty & =  \frac1{N} \sum_{n=1}^N   \E[ \varphi_{n,1}^\top  \varphi_{n,1} ] = \lim_{L\to\infty} \frac1{NL} \sum_{t,n=1}^{L,N}   \varphi_{n,t-1}^\top  \varphi_{n,t-1} \, ,  \\
	\bb_\infty(\cdot, k) & =  \frac1{N} \sum_{n=1}^N   \E[ \varphi_{n,1}^\top c_{n,1}(k)] =  \lim_{L\to\infty} \frac1{NL} \sum_{t,n=1}^{L,N}  \varphi_{n,t-1}^\top c_{n,t}(k), \,\quad  1\leq k\leq K\, , 
\end{aligned}
\end{equation}
where the expectation is with respect to the stationary measure of the Markov chain.

We can obtain the estimator's asymptotical normality using the proof of Theorem \ref{thm:AN_lse} along with the law of large numbers and the central limit theorem for an ergodic trajectory.  
\begin{theorem} \label{thm:AN_lse_1traj} 
Assume that $\bA_\infty $ in \eqref{eq:Abinfty_1traj} is nonsingular. Then, for each $k\in [K]$, the estimator $\widehat \bT_L(\cdot,k)$ in \eqref{eq:Ab_k_1traj} is asymptotically normal, i.e., 
$\sqrt{L} \left( \widehat \bT_L(\cdot,k) -  \bT(\cdot,k) \right) \to  \mathcal{N}(0, \bA_\infty^{-1} \Sigma_k \bA_\infty^{-1} ), 
 $ 
as $L\rightarrow\infty$, where the covariance matrix $\Sigma_k$ is,  with $\overline{\Phi_{n}(k)}: = \lim_{L\to\infty} \frac{1}{L}\sum_{t=1}^L [c_{n,t}(k)\varphi_{n,t-1}]$, 
\begin{equation*}
\begin{aligned}
\Sigma_k =\lim_{L\to\infty}  \frac1{L^2N^2} \sum_{t, t'=1}^L \sum_{n,n'=1}^N \left( c_{n,t}(k)\varphi_{n,t-1}- \overline{\Phi_{n}(k)} \right)\left(  c_{n',t'}(k)\varphi_{n',t'-1}^\top - \overline{\Phi_{n}(k)}^\top  \right) .
\end{aligned}
\end{equation*}
\end{theorem}


\subsection{LSE from ensemble data without trajectory}\label{sec:lse_ens}
Another interesting setting is when the observations are $M_{n,t}$ samples of $X(t)$ for each time $t$, but these samples may come from different trajectories. We call this setting as \emph{ensemble data without trajectory information} and denote the data by 
\[
\textbf{Ensemble Data:}\quad \{X^m_n(t), m=1,\ldots, M_{n,t}\}_{n=1,t=0}^{N,\,\,L}. 
\]
The sample sizes $\{M_{n,t}\}$ do not have to be the same at different times, as long as their minimum $M := \min_{1\leq n\leq N, 1\leq t\leq L} \{ M_{n,t}\}$ is large enough. We will study the error bounds with respect to $M$. 

We estimate $\bT$ by least squares that match the empirical marginal densities of each site. Recall that the marginal density of site $n$ at time $t$ is,  for $1\leq k\leq K$, 
\begin{equation}\label{eq:lse_ens_exp}
p_{n,t} (k) =  \P\{ X_n(t)=k\}  = \E[ \P\{X_n(t)=k | X(t-1) \} ] = \mathbb{E}[ \varphi_{n,t-1}\bT_{\cdot,k}]  = \mathbb{E}[ \varphi_{n,t-1}]\bT_{\cdot,k}.   
\end{equation}
Thus, our estimator is based on empirical approximations of $p_{n,t}(k)$ and $
\mathbb{E}[ \varphi_{n,t-1}]$: 
\begin{equation}\label{eq:ens_empirical}
\begin{aligned}
	 \widehat p_{n,t,M} (k): = \frac{1}{{M_{n,t}}} \sum_{m=1}^{M_{n,t}} \delta_{X^m_n(t)}(k), \quad  
 \widehat \varphi_{n,t-1,M} : =  \frac{1}{M_{n,t-1}} \sum_{m=1}^{M_{n,t-1}}    \varphi_{n,t-1}^m   \in \R^{1\times K} 
\end{aligned}
\end{equation}
for $1\leq k\leq K$ and  $1\leq t\leq L$. 
Note that they are determined by the empirical distributions at each time, and there is no need for sample trajectories.

Our least squares estimator, called \emph{ensemble LSE},  minimizes the discrepancy between the empirical approximations $\widehat p_{n,t+1,M}$ and $\widehat \varphi_{n,t-1,M}$ in \eqref{eq:ens_empirical}:  
\begin{equation}\label{eq:loss_ens}
\begin{aligned}
	\widehat \bT^e_M = \argmin{\bT \in  \R^{K\times K}} \sum_{k=1}^K \sum_{n=1}^N\sum_{t=1}^L 
	\|\widehat p_{n,t,M}(k) - \widehat \varphi_{n,t-1,M} \bT_{\cdot,k}\|^2.
\end{aligned}
\end{equation}
The ensemble LSE is solved by  
\begin{equation}\label{eq:enLSE}
	\begin{aligned}
\widehat {\bT}^e_M( \cdot,k) &= (\bA_{M}^e) ^{\dag}  \bb_{M}^e(\cdot, k), \quad 1\leq k\leq K,  
  \\
\bA^e_{M} & =\frac1{LN}\sum_{t,n=1}^{L,N}   \widehat \varphi_{n,t-1,M}^\top \widehat \varphi_{n,t-1,M}, \quad 
\bb^e_{M}(\cdot,k)  =  \frac1{LN}\sum_{t,n=1}^{L,N} \widehat \varphi_{n,t-1,M}^\top  \widehat p_{n,t,M}(k) \, . 
\end{aligned}
\end{equation}
 Similar to the multi-trajectory LSE in Section \ref{sec:lse_traj}, in practice, we obtain the ensemble LSE by least squares with non-negative constraints, followed by row-normalization.   

This LSE can be viewed as a generalized moment estimator since the entries in the normal matrix and normal vector are approximations of moments. We will show that the estimator is asymptotically normal under a new identifiability condition.

\paragraph{Identifiability in the large sample limit.}
 For $1\leq k\leq K$, denote 
\begin{equation}\label{eq:A_ens}
\begin{aligned}	
	\bA^e_\infty =  \frac1{LN} \sum_{t=1}^L \sum_{n=1}^N   \E[ \varphi_{n,t-1}] ^\top  \E[\varphi_{n,t-1} ], \quad 
	\bb^e_\infty(\cdot, k)  =  \frac1{LN} \sum_{t=1}^L \sum_{n=1}^N   \E[ \varphi_{n,t-1}]^\top p_{n,t}(k) \,. 
\end{aligned}
\end{equation}

\begin{assumption}[Identifiability condition: ensemble data]
\label{assump:ID_ens} 
The distribution of the data satisfies the fact that the matrix $\bA^e_\infty $ in \eqref{eq:A_ens} is non-singular. 
\end{assumption}
 The next lemma shows that the local transition matrix solves the unique solution to the linear system $\bA_\infty^e \bT = \bb_\infty^e$, which follows directly from the definitions of  $\bA_\infty^e$ and $\bb_\infty^e$. 
\begin{lemma}\label{lemma:Ainv_ens} 
Under Assumption {\rm\ref{assump:ID_ens}}, $\bT = (\bA_\infty^e)^{-1}\bb_\infty^e$. 
\end{lemma}

This assumption imposes constraints on both the distribution of the process and the local empirical distributions $\{\varphi_n\}$, which depend on the neighborhood size of the interaction. Three factors can contribute to the identifiability: a non-symmetric initial distribution between sites, a neighborhood that can lead to varying local empirical distributions, and a process that varies in time. For example, $\bA^e_\infty$ can be full rank if $\{ \E[ \varphi_{n,0}]\}_{n=1}^N $ has rank $K$, which relies on a diverse initial distribution and local empirical distribution. Example \ref{exp:full-nbhd} below shows an extreme case that $\{\varphi_n\}$ are the same for all sites, and we rely on the distribution at different times to attain a full-rank norm matrix.   
\begin{example}[Full-network neighborhood] 
\label{exp:full-nbhd}
The local empirical distributions are the same for all sites if the neighborhood is the entire network for each agent. For example, the model in Example {\rm \ref{exp:full-nbhd-N2K2}} has $\varphi_1 = \varphi_2$ for all states $\bx$. 
Thus, for full-network neighborhood, we have $ \E[ \varphi_{n,t}] = \E[ \varphi_{n',t}] $ for any $n\neq n'$. Then, we have $\bA^e_\infty  
 =  \frac1{L} \sum_{t=1}^L\E[ \varphi_{1,t-1}] ^\top  \E[\varphi_{1,t-1} ] $ and it is full rank only if $\{\E[ \varphi_{1,t-1}]\}_{t=1}^L$ has rank $K$. As discussed later, the normal matrix has rank 1 when the process is stationary.  
\end{example}

\paragraph{Asymptotic normality.} We show next that the ensemble LSE is asymptotically normal.

\begin{theorem} \label{thm:AN_ensLSE} 
Under Assumption {\rm\ref{assump:ID_ens}}, for each $ k\in [K]$, the estimator $\widehat \bT_M^e(\cdot,k)$ in \eqref{eq:enLSE} is asymptotically normal;  that is, 
 \begin{equation*}
\sqrt{M} \left( \widehat \bT_M^e(\cdot,k) -  \bT(\cdot,k) \right) \to  \mathcal{N}(0, (\bA_\infty^e)^{-1} \Sigma_k^e(\bA_\infty^e)^{-1} ), 
 \end{equation*}
 where the covariance matrix $\Sigma_k^e$ is, 
\begin{equation}\label{eq:Sigma_k^e}
\begin{aligned}
 \Sigma_k^e =\frac1{L^2N^2} \sum_{n,t=0}^{n,t-1}p_{n,t}(k)^2 \E[ \varphi_{n,t-1}] ^\top  \E[ \varphi_{n,t-1}]) \in \R^{K\times K}. 
\end{aligned}
\end{equation}
\end{theorem}
\begin{proof} 
The proof is based on the asymptotic properties of the empirical approximations of $\widehat p_{n,t,M} (k)$ and $\widehat \varphi_{n,t-1,M} $ defined in \eqref{eq:ens_empirical}. 

First, by the strong Law of Large Numbers, 
$
\widehat \varphi_{n,t-1,M}\xrightarrow{a.s.} \mathbb{E}[ \varphi_{n,t-1}] $
 as $M\to \infty$, for each $n,t$. Also, by the central limit theorem,  
\begin{equation}\label{eq:p_clt}
\sqrt{M} \big(\widehat p_{n,t,M}(k) - p_{n,t}(k) \big) \xrightarrow{d} \mathcal{N}(0,\sigma_{n,t}(k)), 
\end{equation}
for each $n,t$, where the variance $\sigma_{n,t}(k)$ follows from (recall that $\widehat p_{n,t,M} (k)= \frac{1}{{M_{n,t}}} \sum_{m=1}^{M_{n,t}} \delta_{X^m_n(t)}(k)$)  
\begin{equation*}
\begin{aligned}
	\sigma_{n,t}(k) & =  \E [\widehat p_{n,t,M} (k)\widehat p_{n,t,M} (k)] =   \E [  \frac{1}{{M_{n,t}^2}} \sum_{m,m'=1}^{M_{n,t}} \delta_{X^m_n(t)}(k)\delta_{X^{m'}_n(t)}(k)]  \\
	& = \E[ \delta_{X^m_n(t)}(k) ]\E [\delta_{X^{m'}_n(t)}(k)]  = p_{n,t}(k)^2. 
\end{aligned}
\end{equation*}

Next, we study asymptotic properties of the normal matrix and vector in \eqref{eq:enLSE}.  Since $\widehat \varphi_{n,t-1,M}\xrightarrow{a.s.} \mathbb{E}[ \varphi_{n,t-1}]$ for each $n,t$, the normal matrix must also converge a.s., i.e.,
\[
\bA_M^e \xrightarrow{a.s.} \bA^e_\infty = \frac1{LN} \sum_{t=1}^L \sum_{n=1}^N   \E[ \varphi_{n,t-1}] ^\top  \E[ \varphi_{n,t-1}],
\] 
where $\bA^e_\infty$ is defined in \eqref{eq:A_ens}. Meanwhile, by Slutsky's theorem and \eqref{eq:p_clt}, we have
\begin{align*}
\sqrt{M} \widehat\varphi_{n,t-1,M} ^\top [\widehat p_{n,t,M} (k) - p_{n,t} (k)] & \xrightarrow{d}  \mathcal{N}(0, p_{n,t}(k)^2 \E[ \varphi_{n,t-1}] ^\top  \E[ \varphi_{n,t-1}])  
\end{align*}
for each $ n\in [N],0\leq t\leq L-1$. Then, since $\left( \widehat\varphi_{n,t-1,M} ^\top - \E [\widehat\varphi_{n,t-1} ^\top]  \right)  p_{n,t} (k)\xrightarrow{a.s.} 0$, we have
\[
\sqrt{M} \bigg( \widehat\varphi_{n,t-1,M} ^\top \widehat p_{n,t,M} (k) - \E[ \varphi_{n,t-1}^\top] p_{n,t} (k) \bigg ) \xrightarrow{d}  \mathcal{N}(0, p_{n,t}(k)^2 \E[ \varphi_{n,t-1}] ^\top  \E[ \varphi_{n,t-1}]).  
\]
Consequently, 
\[
\sqrt{M} \big( \bb_M^e(\cdot,k) -\bb_\infty^e(\cdot,k) \big )  =  \frac{\sqrt{M}}{LN} \sum_{t=1}^L \sum_{n=1}^N  \bigg( \widehat\varphi_{n,t-1,M} ^\top \widehat p_{n,t,M} (k) - \E[ \varphi_{n,t-1}^\top] p_{n,t} (k) \bigg) \xrightarrow{d} \mathcal{N}(0,\Sigma_k^e), 
\]
with $\Sigma_k^e$ defined in \eqref{eq:Sigma_k^e}.

Then, by Lemma \ref{lemma:Ainv_b_AN} and the invertibility of $\bA^e_\infty$ in Assumption \ref{assump:ID_ens}, we have $\sqrt{M} \bigg( (\bA_M^e)^\dag \bb_M^e(\cdot,k) - (\bA_\infty^e)^{-1} \bb_\infty^e(\cdot,k) \bigg) \xrightarrow{d} \mathcal{N}(0, (\bA_\infty^e)^{-1} \Sigma_k^e (\bA_\infty^e)^{-1})$. This, together with Lemma \ref{lemma:Ainv_ens}, proves the asymptotic normality of the ensemble LSE $\bT_M^e(\cdot,k)$. 
\end{proof}

\subsection{Non-identifiability from stationary distribution}\label{sec:infer_from_inv}
Inference from the stationary distribution is challenging due to the limited information available. It is well known that a stationary distribution of a Markov chain does not determine its transition matrix, i.e., $\bP$ is not determined by $\pi$. Similarly,  the local transition matrix in our model is under-determined by the stationary distribution. Example \ref{ex:1pi-2T} shows that when $(N, K)=(2,2)$, there are two $\bT$'s leading to the same invariant measure. Theorem \ref{thm:nonID-LSE} shows that, in general, the ensemble-LSE is under-determined by the marginal invariant distribution, even though $\bT$ has only $K(K-1)$ unknowns and the invariant measure $\pi$ has $K^N$ entries.

We start with a few basic facts when the process is stationary. By the shift-invariance in Proposition \ref{prop:shift-invariance}, the marginal distributions of all vertices are the same, and so are the expectations of the local empirical distributions, i.e., 
\begin{equation}\label{eq:marginal_local}
	p_{n,t} =p_{1,0}, \quad  \E[ \varphi_{n,t} ] = \E_\pi[\varphi_{1,0}], \quad \forall l\geq 0, \,  n\in [N]. 
\end{equation}  
Thus, the large sample limit of the loss function in \eqref{eq:loss_ens} is 
\begin{equation}\label{eq:lossFn_inv_lse}
\begin{aligned}
	\mathcal{E}_\infty(\bT) & = 
	\sum_{k,n,t=1}^{K,N,L} 
	\| p_{n,t,\infty}(k) -  \varphi_{n,t-1,\infty} \bT_{\cdot,k}\|^2  = \sum_{k=1}^K \big| p_{1,0}(k) -  \E [\varphi_{1,0}] \bT_{\cdot,k}\big|^2 = \| p_{1,0}- \E [\varphi_{1,0}] \bT\|^2_2.
\end{aligned}
\end{equation}

 \begin{theorem}[Non-identifiability from the stationary distributions]
  \label{thm:nonID-LSE}
  Given only the invariant measure, the local transition matrix $\bT$ is under-determined,   i.e., the loss function in \eqref{eq:lossFn_inv_lse} has multiple minimizers,  either when $K>2$ or when $K=2$ with $p_{1,0}= \E[\varphi_{1,0}]$. 
 \end{theorem}

\begin{proof} Note that when given only the stationary measure, the local empirical distributions $\{\varphi_{n}\}$ can only be used via their expectations. Then, by \eqref{eq:marginal_local}, one can only estimate $\bT$ from: (1) the discrepancy between the distributions $p_{1,0}$ and $\E[\varphi_{1,0}]\bT$; and (2) the fact that $p_{1,0} = p_{1,0} \bT_{true}$ by Theorem \ref{thm:1-1PT}. Thus, the identifiability of $\bT$ is equivalent to the uniqueness of the minimizer to the loss function 
 \begin{equation}\label{eq:loss_local}
 	\mathcal{E}_{local}(\bT) =  \| p_{1,0} - p_{1,0} \bT\|_2^2 + \| p_{1,0}- \E [\varphi_{1,0}] \bT\|^2_2. 
 \end{equation} 
	
	Since $\mathcal{E}_{local}(\bT)$ is quadratic in $\bT$, it suffices to study the invertibility of its Hessian 
\begin{equation}\label{eq:A_local_inv}
\begin{aligned}	
	\mathrm{Hess}(\mathcal{E}_{local}) & =  p_{1,0}^\top p_{1,0}+ \E[ \varphi_{1,0}] ^\top  \E[\varphi_{1,0} ]. 
	\end{aligned}
\end{equation}
Here, the Hessian is with respect to $\bT_{\cdot,k}$ and they are the same for all $k\in [K]$. The Hessian matrix has rank $2$ when $p_{1,0}\neq \E[ \varphi_{1,0}]$; and it has rank 1 when $p_{1,0}= \E[ \varphi_{1,0}]$. 
Thus, there are multiple minimizers to $\mathcal{E}_{local}(\bT)$, i.e., $\bT$ is under-determined, either when $K>2$ or when $K=2$ with $p_{1,0}= \E[ \varphi_{1,0}]$. 
\end{proof}

The above non-identifiability is rooted in the limited information available, as only the marginal distributions are used. Loss functions other than the quadratic loss function in \eqref{eq:loss_local}, such as those based on the Kullback-Leibler divergence, total variation, or Wasserstein distances between $p_{1,0}$ and $\E[\varphi_{1,0}] \bT$,  will also have the same issue.

\subsection{Non-asymptotic bounds for the LSEs}
\label{sec:non-asymp}
We establish non-asymptotic bounds for the multi-trajectory LSE $\widehat \bT_M$ in \eqref{eq:traj-LSE} and the ensemble LSE $\widehat \bT_M^e$ in \eqref{eq:enLSE}. Roughly speaking, for small $\epsilon$ and $\delta$, both estimators are $\epsilon$-close to the true local transition matrix with a probability of at least $1-\delta$ when the sample size is of order $M = O(\frac{K^2}{\epsilon^2 \lambda_{min}(\bA_\infty)^2}\ln \frac{K}{\delta})$, but the constant for $\widehat \bT_M^e$ is larger
due to the absence of trajectory information. We postpone the proof to Section \ref{sec-append-nonasym} in the appendix. 
\begin{theorem}\label{thm:concentrationT}
Let $\bT$ be the true local transition matrix. For any $\epsilon, \delta\in (0,1)$, let  $\alpha = \frac{\epsilon}{4} \lambda_{min}(\bA_\infty)$ and $s=\frac{1}{2}  \lambda_{min}(\bA_\infty)  \min\big\{ 1, \frac{\epsilon}{2\|\bT\|_F}\big \} $. The following non-asymptotic bounds hold. 
\begin{itemize}
\item[(a)] 	Under Assumption {\rm\ref{assump:Ainvertible}}, the multi-trajectory LSE $\widehat \bT_M$ in \eqref{eq:traj-LSE} satisfies 
	\begin{equation}\label{eq:non-aysmp}
	\P\{\|\widehat \bT_M - \bT \|_F > \epsilon \}< \delta
	\end{equation}
	if the sample size satisfies  
	$ M > M_{\epsilon, \delta}:= \max\big \{ \frac{24 K^2+4\alpha K}{3\alpha^2}\ln \frac{6K^2}{\delta}, \frac{6+2s}{3s^2}\ln \frac{6K}{\delta} \big \}.  	
	$
 \item[(b)] Under Assumption {\rm\ref{lemma:Ainv_ens}} and assuming  $M_{n,t}\equiv M$, the ensemble LSE $\widehat \bT_M^e$ in \eqref{eq:enLSE} satisfies 
	\begin{equation}\label{eq:non-aysmp_Ens}
	\P\{\|\widehat \bT_M^e - \bT \|_F > \epsilon \}< \delta
	\end{equation}
	if $M$ satisfies  
	 $M > M_{\epsilon, \delta}^e:=  \frac{ 96 K^2+16\alpha K}{3\alpha^2}\ln \frac{12NLK}{\delta}. 	
	$
\end{itemize} 
\end{theorem}

\subsection{Numerical examples}\label{sec:num}

Numerical tests demonstrate that the estimators converge at the rate $M^{-1/2}$ as the sample size increases, in agreement with the theory. They also show that the sampling error may lead to estimators missing the periodic property of the local transition matrix and hence the synchronization; thus, additional techniques, such as an application of a threshold or a sparsity condition, are needed to preserve the additional properties of the local transition matrix that can influence important aspects of the dynamics.   

\begin{figure}[htb]
	\centering
	 \includegraphics[width=0.4\textwidth]{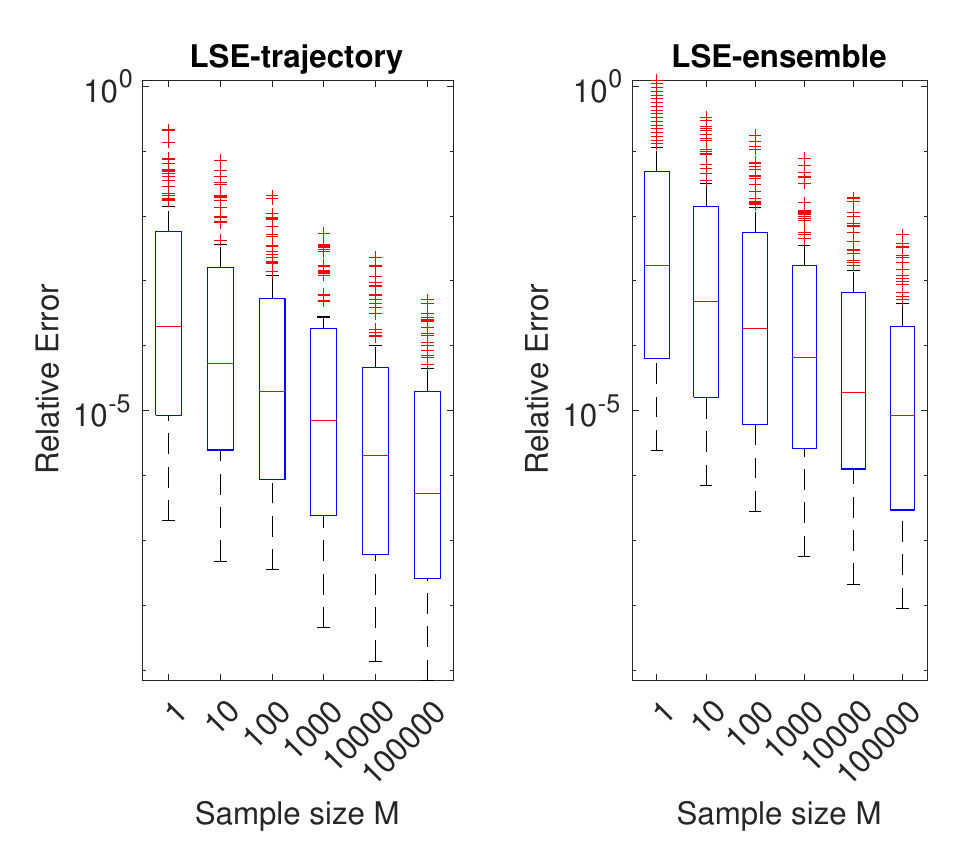}
 	 \includegraphics[width=0.41\textwidth]{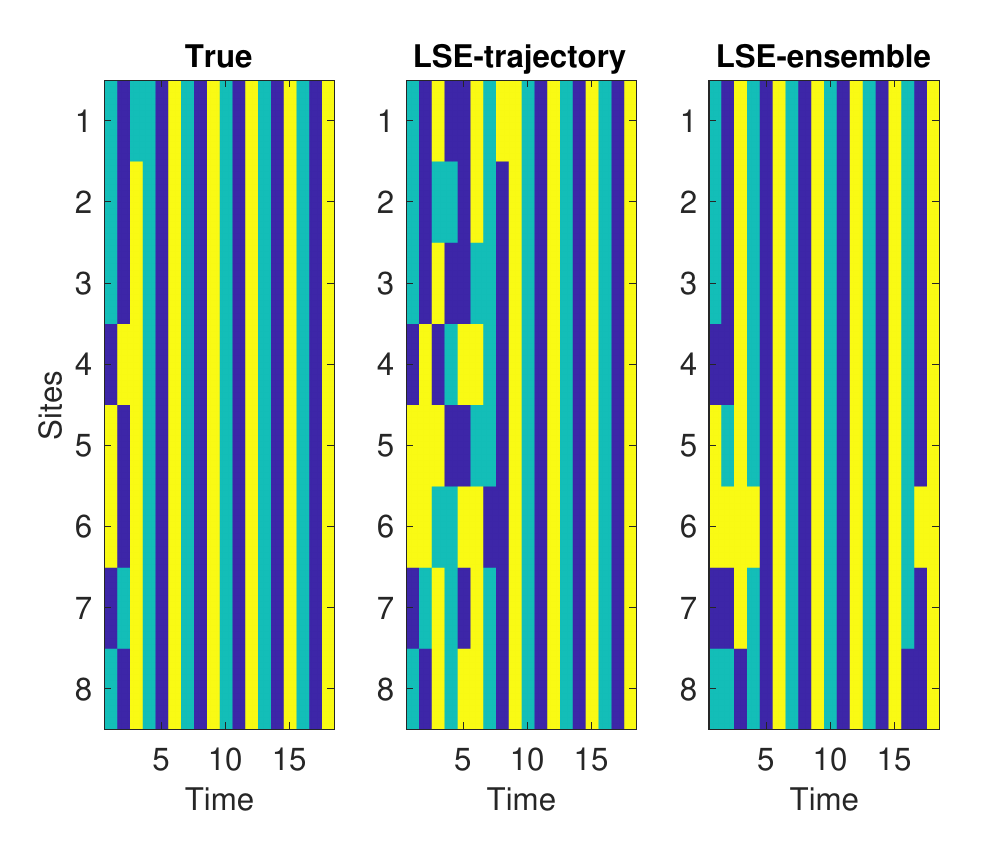}
   \\
        (a) Convergence in sample size \hspace{2cm} (b) Synchronization Prediction\\
	\caption{ \textbf{(a)}: Box plot of relative errors of the LSE estimators in 100 simulations for a system with $(N, K,n_v)=(8,3,3)$. The estimators converge at the same rate; however, the multi-trajectory LSE is significantly more accurate than the ensemble LSE.  \textbf{(b)}: A prediction of synchronization for Example \ref{exp:T-permutation} with $\bT$ estimated from $M=10^3$ trajectories with length $L=100$. The sampling error in LSE-ensemble leads to a system without synchronizations. Here the colors represent the alphabet $\calA= [K] = \{\text{Yellow, Turquoise, and Blue}\}$ with $K=3$.} 
	\label{Fig:lse_conv}
\end{figure}

Figure \ref{Fig:lse_conv}(a) examines the convergence in sample size for the multi-trajectory LSE in Section \ref{sec:lse_traj} and the ensemble LSE in Section \ref{sec:lse_ens}. That is, these estimators are obtained by first solving the normal equations by least squares with non-negative constraints and then row-normalizing the resulting solutions. Here we consider a system with $(N, K,n_v)=(8,3,3)$. The figure shows the box plots of the relative errors of the estimators in 100 independent simulations with increasing sample size. Here we consider a randomly generated matrix $$\bT =  \begin{bmatrix}   0.4719 &    0.0315 &   0.4966 \\
    0.1385  &  0.6118 &   0.2497 \\
    0.2895  &  0.4999   & 0.2107\\
    \end{bmatrix},$$ and use it to generate $5\times 10^5$ sample trajectories with $L=100$. Then, we randomly draw $M$ samples from them $100$ times. 
    
    The results show that the estimators converge as the sample size increases at the rate $M^{-1/2}$, in agreement with Theorems \ref{thm:AN_lse} and \ref{thm:AN_ensLSE}. Additionally, the multi-trajectory LSE is more accurate than the ensemble LSE. As both estimators use the same dataset each time, the increase in accuracy stems from the additional trajectory information.

Figure \ref{Fig:lse_conv}(b) tests the effects of the sampling error in predicting the synchronization. Here the true and estimated local transition matrices are: 
 $$\bT= \begin{bmatrix} 
   0 & 1 & 0 \\
   0 & 0 & 1  \\
   1 &   0 & 0 
\end{bmatrix},  \quad 
\widehat{\bT}_M = \begin{bmatrix} 
0.0000 & 1.0000 & 0.0000 \\
0.0000 & 0.0000 & 1.0000 \\
0.9962 & 0.0038 & 0.0000 
\end{bmatrix}, \quad  
\widehat{\bT}^e_M = \begin{bmatrix} 
 0.0079 & 0.9921 & 0.0000 \\
 0.0182 & 0.0000 & 0.9818 \\
 0.9926 & 0.0074 & 0.0000 
\end{bmatrix}. 
$$
The estimators are inferred from $M=10^2$ sample trajectories with $L=10$. Due to sampling error, neither estimator is periodic; thus, their systems will not synchronize. Figure \ref{Fig:lse_conv}(b) shows that the more accurate multi-trajectory LSE leads to numerical synchronization, while the ensemble LSE cannot maintain the synchronized motion due to the large estimation error. In practice, when the system is known to synchronize, we can apply thresholding or sparsification techniques to preserve the additional properties of the local transition matrix and achieve synchronization.

\section{Future work}
Many venues remain to be explored beyond the scope of this work. 

 The first venue is to study the PCAs on general finite graphs. One may consider general graphs, such as vertex-transitive finite graphs, graphs with a non-negative real-valued weight matrix for edges, or high-dimensional lattices.
The dynamical properties, such as synchronization and ergodicity, as well as the inference of the local transition matrix, can be studied in a similar manner. Additionally, it is of great interest to jointly infer the local transition matrix and the weight matrix of the graph, as studied in \cite{LWLM24} for interacting particle systems on graphs. 

Another venue is to study the PCAs on infinite graphs. Concerning the dynamical properties, one can study the ergodicity and the critical phenomena building on the methods in \cite{toom1994critical,lebowitz1990statistical,casse2023ergodicity,berard2023coupling,follmer2001convergence}. Concerning the inference of the local transition matrix, one may consider the asymptotic and non-asymptotic properties of the estimator when the data is a single trajectory with $N\to \infty$, for which \cite{della2023lan} has established similar results for interacting particle systems and  \cite{bayraktar2024nonparametric} considered this problem for graphon particle systems; see, e.g., \cite{MR4663492}. An interesting parameter to estimate in a similar context would be the size $n_v$ of each neighborhood.

\appendix
\section{Preliminaries and proofs}

\subsection{Properties of Markov chains}\label{sec:prelimMC}
Suppose $X(\cdot)$ is a finite-state Markov chain with an irreducible transition matrix $\bP$. We recall the following general results for Markov chains, and they are used in the proof of Theorem \ref{thm:ergodic}. 
\begin{proposition}
    \label{prop:MC-property}
    \begin{itemize}
        \item It has a unique stationary distribution $\pi$; see {\rm\cite[Theorem 1.7]{durrett2019probability}}. 
        \item All states are positive recurrent; see {\rm\cite[Theorem 1.30]{durrett2019probability}}. 
        \item Suppose $X(\cdot)$ is aperiodic. Then $\lim_{t \to \infty} \prob{X(t) = \cdot} = \pi$; see {\rm\cite[Theorem 1.19]{durrett2019probability}}. 
    \end{itemize}
\end{proposition}

Regarding the exponential convergence to the stationary distribution, we have the following result taken from \cite[Theorem 1.3]{kulik2015introduction}. 

\begin{proposition}
    \label{prop:MC-geometric-cvg}
    Suppose there is some $\rho < 1$ such that
    \begin{equation}
        \label{eq:uniform-ergodic-condition}
        \|\bP(x,\cdot)-\bP(x',\cdot)\|_{TV} \le 2\rho, \quad \forall\,x,x'.
    \end{equation}
    Then
    $\|\bP_t(x,\cdot)-\bP_t(x',\cdot)\|_{TV} \le 2\rho^t$ and 
    $\|\bP_t(x,\cdot)-\pi\|_{TV} \le 2\rho^t$,  $\forall\,t \ge 1, \: x,x'$.
\end{proposition}

\begin{proof}[Proof of Theorem \ref{thm:ergodic}]
    By Proposition \ref{prop:irreduciblePT},  $X(\cdot)$ is irreducible and aperiodic.
    Then, it remains to prove \eqref{eq:geometic-cvg}, as the rest of the statements are classic results (see, e.g., Proposition \ref{prop:MC-property}).

    Since $X(\cdot)$ is irreducible and aperiodic, there exists some $t_0 \in \mathbb{N}$ such that
        $\bP^{t_0}(x,y)>0, \, \forall\, x,y \in [K]^N$.
    Therefore, $|\bP^{t_0}(x,y)-\bP^{t_0}(x',y)|< \max\{\bP^{t_0}(x,y), \bP^{t_0}(x',y)\}<\bP^{t_0}(x,y)+\bP^{t_0}(x',y)$ for any $x,x', y\in [K]^N$. Hence, $$\rho_0:= \frac{1}{2}\max_{x,x' \in [K]^N} \sum_{y\in [K]^N} \max\{\bP^{t_0}(x,y), \bP^{t_0}(x',y)\}  < \frac{1}{2} \max_{x,x' \in [K]^N} \sum_{y\in [K]^N} [ \bP^{t_0}(x,y)+\bP^{t_0}(x',y)]=  1$$ and 
    \begin{align*}
        \|\bP^{t_0}(x,\cdot)-\bP^{t_0}(x',\cdot)\|_{TV} & = \sum_{y\in [K]^N} |\bP^{t_0}(x,y)-\bP^{t_0}(x',y)|
         <  \sum_{y\in [K]^N} \max\{\bP^{t_0}(x,y), \bP^{t_0}(x',y)\} \leq   2\rho_0\, 
    \end{align*}
    for all $x,x' \in [K]^N$.
    It then follows from Proposition \ref{prop:MC-geometric-cvg} that
    $
        \|\bP^{kt_0}(x,\cdot)-\pi\|_{TV} \le 2\rho_0^{k}, \quad \forall\,k \in \mathbb{N},\: x \in [K]^N
    $.
    For $t \ge t_0$, writing $t = \lfloor \frac{t}{t_0} \rfloor t_0 + (t-\lfloor \frac{t}{t_0} \rfloor t_0)$, we have
    \[
        \|\bP^t(x,\cdot)-\pi\|_{TV} \le \|\bP^{\lfloor \frac{t}{t_0} \rfloor t_0}(x,\cdot)-\pi\|_{TV} \le 2\rho_0^{\lfloor \frac{t}{t_0} \rfloor} = \frac{2}{\rho_0} \rho_0^{\lfloor \frac{t}{t_0} \rfloor + 1}\le \frac{2}{\rho_0} (\rho_0^{1/t_0})^t, \quad \forall\, x \in [K]^N.
    \]
    This gives \eqref{eq:geometic-cvg} with $C= \frac{2}{\rho_0} \in (0,\infty)$ and $\rho=\rho_0^{1/t_0} \in (0,1)$.
\end{proof}

\subsection{Concentration inequalities}
\begin{theorem}[Bernstein's Inequality]\label{thm:Berstein} {\rm (see, e.g., \cite[Theorem 2.8.4]{vershynin2018})}
Let $X_1, \dots, X_M$ be independent zero-mean random variables. Suppose that $|X_i| \leq c$ almost surely for all $i$. Then for all positive $t$, 
$\mathbb{P}\left(|\sum_{i=1}^{M}X_i| \geq t \right) \leq 2\exp\left(-\frac{ t^2/2}{\operatorname{Var}(\sum_i X_i) + \frac{1}{3}c t}\right)
$. 
In particular, when $\{X_i\}$ are iid., we have $\mathbb{P}\left(|\frac{1}{M}\sum_{i=1}^{M}X_i| \geq t \right) \leq 2\exp\left(-\frac{ Mt^2/2}{\operatorname{Var}( X_1) + \frac{1}{3}c t}\right)$.   
\end{theorem}

\begin{theorem}[Matrix Bernstein's inequality]{\rm(\cite[Theorem 5.4.1]{vershynin2018}  or \cite[Theorem 6.1.1]{Tropp2015})}
\label{thm:matBerstein}
	Let $\{X_i\}_{i=1}^M\subset \R^{n\times n}$ be independent mean zero symmetric random matrices such that $\|X_i\|_{op}\leq c$ almost surely for all $i$. Then, for every $t\geq 0$, we have 
	$
	\P(\|\sum_{i=1}^M X_i\|_{op}\geq t ) \leq 2n \exp\left( -\frac{t^2/2}{\sigma^2+c t/3}\right), 
	$
	where $\sigma^2= \|\sum_{i=1}^M \E[X_i^2]\|_{op}$. Additionally, when $\{X_i\}$ are identically distributed, we have 
	\[
	\P(\|\frac{1}{M}\sum_{i=1}^M X_i\|_{op}\geq t ) \leq 2n \exp\left( -\frac{M t^2/2}{\|\E[X_1^2]\|_{op}+c t/3}\right). 
	\]  
	
\end{theorem}

\subsection{Proof of Theorem \ref{thm:Lip_P_T}}\label{sec:Proof-Lip_P_T} 
The proof of Theorem \ref{thm:Lip_P_T} is based on the following result from \cite{seneta1988perturbation},  
 which bounds the impact on the invariant measure when perturbing a transition matrix.  
We refer to the general study on the perturbation of Markov chains in \cite{schweitzer1968perturbation,haviv1984perturbation,cho2001comparison}. 
\begin{lemma}
    \label{lem:perturbation-MC}
    For two finite irreducible transition matrices $\bP_1$ and $\bP_2$ with stationary distributions $\pi_1$ and $\pi_2$,
    $\|\pi_1-\pi_2\|_1 \le \frac{1}{1-\tau(\bP_1)} \sum_x \|\bP_1(x,\cdot) - \bP_2(x,\cdot)\|_1.$
\end{lemma}

\begin{proof}[Proof of Theorem \ref{thm:Lip_P_T}] 
    For $x=(x_1,\dotsc,x_N)\in [K]^N$ and $y=(y_1,\dotsc,y_N) \in [K]^N$, writing
    \[
    z_n \equiv z_n(x,y) = \frac{1}{|V_n|} \sum_{n'\in V_n} \bT^{(1)}_{x_{n'},y_n}, \qquad w_n \equiv w_n(x,y) = \frac{1}{|V_n|} \sum_{n'\in V_n} \bT^{(2)}_{x_{n'},y_n},
    \]
    we have 
    $
    \bP_1(x,y) = \prod_{n=1}^N z_n$, $\bP_2(x,y) = \prod_{n=1}^N w_n
    $.
    By adding and subtracting terms,
    \begin{align*}
        |\bP_1(x,y)-\bP_2(x,y)| & = \left|\sum_{n=1}^N w_1 \dotsm w_{n-1} (z_n-w_n) z_{n+1} \dotsm z_N\right|
         \le \sum_{n=1}^N w_1 \dotsm w_{n-1} |z_n-w_n| z_{n+1} \dotsm z_N.
    \end{align*}
    Noting that $z_n,w_n$ depend on $y$ only through $y_n$ and $\bT^{(1)}$ and $\bT^{(2)}$ are stochastic matrices, we have
    $
    \sum_{y_n\in [K]} z_n = 1$, $\sum_{y_n\in [K]} w_n = 1
    $.
    Therefore
    \begin{align*}
        \|\bP_1-\bP_2\|_1 & = \sum_{x,y} |\bP_1(x,y)-\bP_2(x,y)| \le \sum_{n=1}^N \sum_{x,y} w_1 \dotsm w_{n-1} |z_n-w_n| z_{n+1} \dotsm z_N \\
        & = \sum_{n=1}^N \sum_{x} \left(\sum_{y_1} w_1\right) \dotsm \left(\sum_{y_{n-1}}w_{n-1}\right) \left(\sum_{y_n}|z_n-w_n|\right) \left(\sum_{y_{n+1}}z_{n+1}\right) \dotsm \left(\sum_{y_N}z_N\right) \\
        & = \sum_{n=1}^N \sum_{x,y_n} |z_n-w_n| \le \sum_{n=1}^N \sum_{x,y_n} \frac{1}{|V_n|} \sum_{n'\in V_n} |\bT^{(1)}_{x_{n'},y_n} - \bT^{(2)}_{x_{n'},y_n}|. 
    \end{align*}
    Noting that by symmetry, $\sum_{x} \frac{1}{|V_n|} \sum_{n'\in V_n} |\bT^{(1)}_{x_{n'},y_n} - \bT^{(2)}_{x_{n'},y_n}| = K^{N-1}\sum_{x_n} |\bT^{(1)}_{x_{n},y_n} - \bT^{(2)}_{x_{n},y_n}|$, we can write the last term as
    \[
    \sum_{n=1}^N K^{N-1} \sum_{x_n,y_n} |\bT^{(1)}_{x_n,y_n} - \bT^{(2)}_{x_n,y_n}| = \sum_{n=1}^N K^{N-1} \|\bT^{(1)}-\bT^{(2)}\|_1 = N K^{N-1} \|\bT^{(1)}-\bT^{(2)}\|_1.
    \]
    Alternatively, note that
    \[
    \sum_{y_n} |z_n-w_n| = \sum_{y_n} |\varphi_n(\cdot) (\bT^{(1)}_{\cdot,y_n}-\bT^{(2)}_{\cdot,y_n})| \le \|\bT^{(1)}-\bT^{(2)}\|_{1,1}.
    \]
    So, we also obtain 
   $    \|\bP_1-\bP_2\|_1 \le NK^N \|\bT^{(1)}-\bT^{(2)}\|_{1,1}$ and \eqref {eq:DeltaP-L1}.

   To prove \eqref{eq:DeltaP-L2}, using Cauchy--Schwarz inequality, we obtain
    \begin{align*}
        & |\bP_1(x,y)-\bP_2(x,y)|^2  = \big|\sum_{n=1}^N w_1 \dotsm w_{n-1} (z_n-w_n) z_{n+1} \dotsm z_N\big|^2\\
        & \le N\sum_{n=1}^N |w_1 \dotsm w_{n-1} (z_n-w_n) z_{n+1} \dotsm z_N|^2 = N\sum_{n=1}^N w_1^2 \dotsm w_{n-1}^2 |z_n-w_n|^2 z_{n+1}^2 \dotsm z_N^2.
    \end{align*}
    Since $ \sum_{y_n} z_n^2 =\sum_{y_n} \big |  \frac{1}{|V_n|} \sum_{n'\in V_n} \bT^{(1)}_{x_{n'},y_n} \big |^2 \leq  \sum_{y_n} \frac{1}{|V_n|} \sum_{n'\in V_n} \big | \bT^{(1)}_{x_{n'},y_n} \big |^2$ and similarly for $w_n$, i.e., 
    $
    \sum_{y_n} z_n^2 \le C_K$ and $\sum_{y_n} w_n^2 \le C_K,
    $
    we have
    \begin{align*}
        \|\bP_1-\bP_2\|_2^2 & = \sum_{x,y} |\bP_1(x,y)-\bP_2(x,y)|^2 \le N\sum_{n=1}^N \sum_{x,y} w_1^2 \dotsm w_{n-1}^2 |z_n-w_n|^2 z_{n+1}^2 \dotsm z_N^2 \\
        & \le N C_K^{N-1} \sum_{n=1}^N \sum_{x,y_n} |z_n-w_n|^2 \le N C_K^{N-1} \sum_{n=1}^N \sum_{x,y_n} \frac{1}{|V_n|} \sum_{n'\in V_n} |\bT^{(1)}_{x_{n'},y_n} - \bT^{(2)}_{x_{n'},y_n}|^2. 
    \end{align*}
Then,  we obtain \eqref{eq:DeltaP-L2} by using the symmetry to write the last term as 
    \begin{align*}
        N C_K^{N-1} \sum_{n=1}^N K^{N-1} \sum_{x_n,y_n} |\bT^{(1)}_{x_n,y_n} - \bT^{(2)}_{x_n,y_n}|^2 & = N C_K^{N-1} K^{N-1} \sum_{n=1}^N \|\bT^{(1)}-\bT^{(2)}\|_2^2 \\
        & = N^2 (KC_K)^{N-1} \|\bT^{(1)}-\bT^{(2)}\|_2^2.
    \end{align*}

    Finally we estimate $\|\pi_1 - \pi_2\|_1$.
    Note that $\bP_1(x,\cdot)$ may be singular with respect to $\bP_1(x',\cdot)$ for some $x$ and $x'$, resulting in
    $\tau(\bP_1)=1$ which prevents applying Lemma \ref{lem:perturbation-MC} directly.
    However, $\bP_1$ is irreducible and aperiodic, and hence there exists some $\ell_0 \in \mathbb{N}$ such that $\bP_1^{\ell_0}(x,y) > 0$ for all $x,y$ and hence $\tau(\bP_1^{\ell_0})<1$.
    Since $\pi_1$ and $\pi_2$ are also stationary distributions of $\bP_1^{\ell_0}$ and $\bP_2^{\ell_0}$, Lemma \ref{lem:perturbation-MC} implies
    \begin{equation}
    \label{eq:Lipschitz_pf_1} 
    \|\pi_1 - \pi_2\|_{1} \le \frac{1}{1-\tau(\bP_1^{\ell_0})} \|\bP_1^{\ell_0} - \bP_2^{\ell_0}\|_{1}.
    \end{equation}
    By adding and subtracting terms,
    \[
    \|\bP_1^{\ell_0} - \bP_2^{\ell_0}\|_{1} \le \sum_{\ell=1}^{\ell_0} \|\bP_2^{\ell-1}\bP_1^{\ell_0-\ell+1} - \bP_2^{\ell}\bP_1^{\ell_0-\ell}\|_{1} = \sum_{\ell=1}^{\ell_0} \|\bP_2^{\ell-1}(\bP_1-\bP_2)\bP_1^{\ell_0-\ell}\|_{1}.
    \]
    Note that for an $n \times n$ matrix $\bA$ and an $n \times n$ stochastic matrix $\bB$, 
    \begin{align*}
    \|\bA\bB\|_1 &= \sum_{i,j} |(\bA\bB)_{i,j}| \le \sum_{i,j,k} |\bA_{ik}||\bB_{kj}| = \sum_{i,k} |\bA_{ik}| = \|\bA\|_1, \\
    \|\bB\bA\|_1 &= \sum_{i,j} |(\bB\bA)_{i,j}| \le \sum_{i,j,k} |\bB_{ik}||\bA_{kj}| \le \sum_{j,k} |\bA_{kj}| \max_{k'} \sum_{i} |\bB_{ik'}| \le n \|\bA\|_1.    	
    \end{align*}   
    Applying the above two estimates with $n=K^N$, we have
    \[
    \|\bP_1^{\ell_0} - \bP_2^{\ell_0}\|_{1} \le \sum_{\ell=1}^{\ell_0} \|\bP_2^{\ell-1}(\bP_1-\bP_2)\|_{1} \le (1+(\ell_0-1)K^N)\|\bP_1-\bP_2\|_{1},
    \]
    where $\bB=\bP_1^{\ell_0-\ell}$ for the first inequality and $\bB=\bP_2^{\ell-1}$ for the second inequality.
    Combining this with \eqref{eq:Lipschitz_pf_1}  gives the desired result.
\end{proof}

\subsection{Proof for the non-asymptotic bounds}\label{sec-append-nonasym}

The proof is based on the concentration bounds for the normal matrices and vectors in the next lemma. These bounds highlight that the trajectory-based normal matrice and vector approach their large sample limits faster than those without using trajectory information. 
\begin{lemma}[Concentration for normal matrices and vectors]
\label{eq:concentrationAb} 
 For any $s>0$, the following concentration bounds hold for the normal matrix $\bA_M$ and vector $\bb_M$ in for the multi-trajectory LSE in \eqref{eq:traj-LSE}, and $\bA_M^e$ and $\bb_M^e$ for the ensemble LSE \eqref{eq:enLSE}. 
 \begin{itemize}
 \item[(a)] Under Assumption {\rm\ref{assump:Ainvertible}}, we have
 	\begin{equation*}
\begin{aligned}
		  \prob{ \|\bA_\infty- \bA_M\|_{op} >  s} & < 2K \exp\left(-\frac{Ms^2/2}{1+ s/3} \right), \\
		  \prob{ \|\bb_\infty- \bb_M\|_{F} >  s} & < 2K^2 \exp\left(-\frac{Ms^2}{8 K^2+ 4sK/3} \right). \\
\end{aligned}
\end{equation*}
\item[(b)] Under Assumption {\rm\ref{assump:ID_ens}}, we have
 	\begin{equation*}
\begin{aligned}
		  \prob{ \|\bA_\infty^e- \bA_M^e\|_{op} >  s} & < 2NLK \exp\left(-\frac{Ms^2}{ 32K^2+ 16Ks/3} \right), \\
		  \prob{ \|\bb_\infty^e- \bb_M^e\|_{F} >  s} & < 4NLK \exp\left(-\frac{Ms^2}{32 K^2+ 16Ks/3} \right). \\
\end{aligned}
\end{equation*}
 \end{itemize}
\end{lemma}

\begin{proof} These bounds follow from applying Bernstein's inequalities. 

\textbf{Part (a)}. First, note that 
 $ \bA_M-\bA_\infty= \frac{1}{M} \sum_{m=1}^M \left( \bA_{L,N}^m - \E[\bA_{L,N}^m] \right)$, where $\{\bA_{L,N}^m \}_{m=1}^M$, 
 defined in \eqref{eq:Ab-m}, is a sequence of symmetric identically distributed random matrices with  
\begin{align*}
\sigma_1^2& = \| \E[\left( \bA_{L,N}^m- \E[\bA_{L,N}^m]\right)^2] \|_{op} 
               =  \left\| \E[ (\bA_{L,N}^m)^2] - \E[\bA_{L,N}^m]^2  \right\|_{op}^2 \\ 
          & \leq  \E \left\| \bA_{L,N}^m \right\|_{op}^2+ \left\|\E[\bA_{L,N}^m]\right\|_{op}^2  \leq 2 \E\|\bA_{L,N}^m\|_{op}^2. 	
\end{align*}
Meanwhile, since $\bA_{L,N}^m$ is symmetric, 
\begin{align*}
\|\bA_{L,N}^m\|_{op} & = \sup_{u\in \R^{K\times 1}, \|u\|=1} u^\top \bA_{L,N}^m u =  \sup_{u\in \R^{K}, \|u\|=1} \frac1{LN}\sum_{t,n=1}^{L,N}  u^\top ( \varphi_{n,t-1}^{m})^\top   \varphi_{n,t-1}^{m} u \, \\
& = \sup_{u\in \R^{K}, \|u\|=1} \frac1{LN}\sum_{t,n=1}^{L,N}  |\varphi_{n,t-1}^{m} u|^2 \leq 1,  	
\end{align*}
where the inequality follows from $ |\varphi_{n,t-1}^{m} u|^2\leq \|u\|_2  $ since each entry of $\varphi_{n,t-1}^{m}$ is in $[0,1]$. Consequently, $\sigma_1^2\leq 2$. Applying the matrix Bernstein's inequality (see Theorem \ref{thm:matBerstein}), we obtain the bound for $\prob{ \|\bA_\infty- \bA_M\|_{op} >  s}$. 

Next, recall that by the definition of $\bb_M$ in \eqref{eq:traj-LSE}, we have 
\[
[\bb_\infty- \bb_M](k',k) = \frac{1}{M}\sum_{m=1}^M \xi^m_{k',k}\, , \quad \xi^m_{k',k}:= \left( \E [\bb_{L,N}^m(k',k)] - \bb_{L,N}^m(k',k)\right)
\]
with $\bb_{L,N}^m(k',k) =\frac1{LN}\sum_{t,n}^{L,N} c_{n,t}^m(k)   \varphi_{n,t-1}^{m,k'}$. Using the fact that $c_{n,t}^m(k) \in \{0,1\}$ and $ \varphi_{n,t-1}^{m,k'}\in [0,1]$, we have  $|\xi^m_{k',k}|\leq 2$ and $\E[|\xi^m_{k',k}|^2]\leq 4$. Thus, Bernstein's inequality (see Theorem \ref{thm:Berstein}) implies 
\[
\prob{ |[ \bb_\infty- \bb_M](k',k)| > \frac{s}{K} } < 2\exp\left( -\frac{Ms^2/(2K^2)}{4+2s/(3K)} \right) = 2\exp\left( -\frac{Ms^2}{8K^2+4sK/3)} \right). 
\]
Meanwhile, note that if $\|\bb_\infty- \bb_M\|_F^2 = \sum_{k',k=1}^{K,K} | \bb_\infty- \bb_M](k',k)|^2 >s^2$ holds true, then there exists $k',k$ such that $| \bb_\infty- \bb_M](k',k)| > \frac{s}{K} $. Hence, 
\begin{align*}
\prob{ \|\bb_\infty- \bb_M\|_{F} >  s} & \leq \prob{ \bigcup_{k',k}\{| [\bb_\infty- \bb_M](k',k)| > \frac{s}{K}\} }  \\
& \leq \sum_{k',k} \prob{ | [\bb_\infty- \bb_M](k',k)| > \frac{s}{K} } \leq  2K^2\exp\left( -\frac{Ms^2}{8K^2+4sK/3)} \right). 
\end{align*}

\textbf{Part (b)}.  The normal matrix $\bA_M^e$ and vector $\bb_M^e$ in \eqref{eq:enLSE} require additional treatments since they are the averages of products of sample means, and this prevents directly applying concentration inequalities. Recall that $\bA_M^e= \frac1{LN}\sum_{t,n=1}^{L, N}   \widehat \varphi_{n,t-1, M}^\top \widehat \varphi_{n,t-1, M} \in \R^{K\times K}$ and $\bb_M^e= \frac1{LN}\sum_{t,n=1}^{L, N} \widehat \varphi_{n,t-1, M}^\top  \widehat p_{n,t, M}$,  where  $ \widehat \varphi_{n,t,M} =  \frac{1}{M} \sum_{m=1}^{M}    \varphi_{n,t}^m $ and  $\widehat p_{n,t,M}  = \frac{1}{M} \sum_{m=1}^{M} p_{n,t}^m$ are sample averages defined in \eqref{eq:ens_empirical}. Here, we denote $p_{n,t}^m\in \{0,1\}^K$ with $p_{n,t}^m (k): = \delta_{X^m_n(t)}(k)$  for $1\leq k\leq K$.   

To apply the concentration inequality, we use the upper bounds of $\widehat \varphi_{n,t-1, M}$ or $\widehat p_{n,t, M}$ to remove one of the terms in the product so that we can apply concentration inequality to the remaining sample mean term. This will lead to a multiplicative factor $NL$ in the upper bounds for the probabilities.  

First, we show the concentration bounds for the following centered random vectors:  \vspace{-2mm}
\begin{equation*}
\widetilde \varphi_{n,t,M} := \widehat \varphi_{n,t,M}- \E[\widehat \varphi_{n,t,M}] = \frac{1}{M} \sum_{m=1}^{M}  \widetilde \varphi_{n,t}^m,\quad  \widetilde p_{n,t,M} := \widehat p_{n,t,M}- \E[\widehat p_{n,t,M}] = \frac{1}{M} \sum_{m=1}^{M}  \widetilde p_{n,t}^m, 
\end{equation*} 
where $\widetilde \varphi_{n,t}^m:= \varphi_{n,t}^m - \E[\varphi_{n,t}^m ]$ and $\widetilde p_{n,t}^m:= p_{n,t}^m - \E[p_{n,t}^m ]$, and show that: \vspace{-2mm}
\begin{equation}\label{eq:phi_p_tilde}
\begin{aligned}
 \| \widetilde \varphi_{n,t,M}\|_{2} \leq \sum_{k=1}^K |\widetilde \varphi_{n,t,M}(k)|\leq 2, \quad  \quad  \| \widetilde p_{n,t,M}\|_{2} \leq   \sum_{k=1}^K |\widetilde p_{n,t,M}(k)|\leq 2. 
\end{aligned}
\end{equation}

The arguments for $\widetilde \varphi_{n,t,M} $ and  $\widetilde p_{n,t,M}$ are identical, so we only consider $\widetilde \varphi_{n,t,M} $.

 Since $ \varphi_{n,t}^m$ is a probability distribution, i.e., its entries are non-negative and $\sum_{k=1}^K \varphi_{n,t}^m(k)=1$, we have $ |\widetilde \varphi_{n,t}^m (k)|= | \varphi_{n,t}^m(k) - \E[ \varphi_{n,t}^m ](k) | \leq  \varphi_{n,t}^m(k) + \E[ \varphi_{n,t}^m(k) ] \leq 2 $ for each $k$ and $\sum_k |\widetilde \varphi_{n,t}^m (k)| \leq 2$. 
Applying that $\|a\|_{2 }= \left( \sum_{k=1}^K a_k^2 \right)^{1/2}\leq \sum_{k=1}^K |a_k|$, we obtain \eqref{eq:phi_p_tilde} from  
\begin{align*}
	 \| \widetilde \varphi_{n,t,M} \|_2  \leq \sum_{k=1}^K |\widetilde \varphi_{n,t,M}(k)| &= \sum_{k=1}^K \left| \frac{1}{M} \sum_{m=1}^M  \widetilde \varphi_{n,t}^m (k)  \right| \leq \sum_{k=1}^K \frac{1}{M} \sum_{m=1}^M  | \widetilde \varphi_{n,t}^m (k) |   
    \leq  2. 
\end{align*}

Also, for each $k$, applying Bernstein's inequality (Theorem \ref{thm:Berstein}) to $\widetilde \varphi_{n,t,M}(k)=  \frac{1}{M} \sum_{m=1}^{M}  \widetilde \varphi_{n,t}^m(k) $ with $|\widetilde \varphi_{n,t}^m (k)|\leq 2$ and $\mathrm{Var}(\widetilde \varphi_{n,t}^m (k))\leq 4$, we obtain (and similarly for $\widetilde p_{n,t,M}(k)$): 
\begin{equation}\label{eq:phi_p_concentr}
 \prob{ | \widetilde \varphi_{n,t,M}(k)| >\frac{s}{4K}  } \leq  2e^{ -\frac{Ms^2}{32K^2+16sK/3 } }, \quad  \prob{ | \widetilde p_{n,t,M}(k)| >\frac{s}{4K}  } \leq  2e^{ -\frac{Ms^2}{32K^2+16sK/3 } }.   	
\end{equation}

Next, we show the bound for $\| \bA_M^e-\bA_\infty^e\|_{op}$. Using $\widehat \varphi_{n,t-1,M} =\widetilde \varphi_{n,t,M}+ \E[\widehat \varphi_{n,t,M}] $,  we have 
\begin{align*}
	\bA_{M}^e - \bA_\infty^e  & = \frac1{NL}\sum_{t,n=1}^{L,N}  \bigg( \widehat \varphi_{n,t-1,M} ^\top\widehat \varphi_{n,t-1,M}  - \E[\widehat \varphi_{n,t-1,M} ^\top\widehat \varphi_{n,t-1,M}] \bigg) \\
	& =  \frac1{NL}\sum_{t,n=1}^{L,N} \left(  \widetilde \varphi_{n,t-1,M} \widetilde \varphi_{n,t-1,M}^\top + \widetilde \varphi_{n,t-1,M}^\top \E [\widehat \varphi_{n,t-1,M} ] + \E [\widehat \varphi_{n,t-1,M} ]^\top \widetilde \varphi_{n,t-1,M} \right). 
\end{align*}
 Then,  noting that $\| \widetilde \varphi_{n,t-1,M}^\top \E [\widehat \varphi_{n,t-1,M} ]\|_{op} = \|\E [\widehat \varphi_{n,t-1,M} ]^\top \widetilde \varphi_{n,t-1,M}\|_{op}$ and that for any $u,v\in \R^{1\times K}$,  $\| u^\top v\|_{op}= \sup_{c\in \R^{1\times K}, \|c\|_2=1} cu^\top vc^\top \leq \|u\|_2\|v\|_2 $, and , and we have 
\begin{align*}
	\| \bA_{M}^e - \bA_\infty^e  \|_{op} 
	& \leq  \frac1{NL}\sum_{t,n=1}^{L,N} \|  \widetilde \varphi_{n,t-1,M}^\top \widetilde \varphi_{n,t-1,M} \|_{op} +2 \| \widetilde \varphi_{n,t-1,M}^\top \E [\widehat \varphi_{n,t-1,M} ]\|_{op}
	\\
	&\leq  \frac1{NL}\sum_{t,n=1}^{L,N} (  \|\widetilde \varphi_{n,t-1,M} \|_2^2 + 2 \|\widetilde \varphi_{n,t-1,M} \|_2 \|\E[\widetilde \varphi_{n,t-1,M}]\|_2) \\
	& \leq  \frac 4 {NL}\sum_{t,n=1}^{L,N}  \|\widetilde \varphi_{n,t-1,M} \|_2 \leq  \frac 4 {NL}\sum_{t,n=1}^{L,N} \sum_{k=1}^K | \widetilde \varphi_{n,t-1,M}(k)|, 
\end{align*}
where the last two inequalities follow from  \eqref{eq:phi_p_tilde} and $ \| \E[\widehat \varphi_{n,t,M}] \|_2 \leq \sum_k \E[\widehat \varphi_{n,t,M}(k)]  = 1$ for any $n,t,M$. 
Then, \eqref{eq:phi_p_concentr} implies 
\begin{align*}
\prob{ \|\bA_M^e -\bA_\infty^e\|_{op} >  s} & \leq   \prob{ \frac 4 {NL}\sum_{t,n=1}^{L,N} \sum_{k=1}^K | \widetilde \varphi_{n,t-1,M}(k)| >s }  \leq \prob{ \bigcup_{t,n,k=1}^{L,N,K}\{ |\widetilde \varphi_{n,t-1,M}(k)|| >\frac{s}{4K} \}}\\
& \leq \sum_{t,n=1}^{L,N} \sum_{k} \prob{ | \widetilde \varphi_{n,t-1,M}(k)| >\frac{s}{4K}  } \leq  2NLK\exp\left( -\frac{Ms^2}{32K^2+16sK/3)} \right). 
\end{align*}
Similarly, the same bound holds for $\widetilde p_{n,t,M}$. 

Lastly, consider $\|\bb_M^e -\bb_\infty^e\|_{F}$. Note that $ \| \E[\widehat p_{n,t,M}] \|_2 \leq \sum_k \E[\widehat p_{n,t,M}(k)]  = 1$ and $ \|\E[\widehat \varphi_{n,t,M}] \|_2 \leq \sum_k \E[\widehat \varphi_{n,t,M}(k)]  = 1$ for any $(n,t,M)$. Also, $\|\widetilde \varphi_{n,t-1,M} \|_2\|\widetilde p_{n,t,M} \|_2\leq \|\widetilde \varphi_{n,t-1,M} \|_2 + \|\widetilde p_{n,t,M} \|_2 $ by \eqref{eq:phi_p_tilde}. Using $\|u^\top v\|_F = \|u\|_2\|v\|_2 $ for any $u,v\in \R^{1\times K}$ (since $\|u^\top v\|_F^2 = \|u\|_2^2\|v\|_2^2 $),  we have 
\begin{align*}
  &\| \bb_{M}^e - \bb_\infty^e  \|_F 
        \leq  \frac1{NL}\sum_{t,n=1}^{L,N} \|  \widetilde \varphi_{n,t-1,M}^\top \widetilde p_{n,t,M} \|_F + \| \widetilde \varphi_{n,t-1,M}^\top \E [\widehat p_{n,t,M} ]\|_F + \|\E [\widehat \varphi_{n,t-1,M} ]^\top \widetilde p_{n,t,M}\|_F \\
	&=  \frac1{NL}\sum_{t,n=1}^{L,N} (  \|\widetilde \varphi_{n,t-1,M} \|_2\|\widetilde p_{n,t,M} \|_2 +  \|\widetilde \varphi_{n,t-1,M} \|_2 \|\E [\widehat p_{n,t,M}]\|_2 + \|\E[\widehat \varphi_{n,t-1,M}]\|_2\|\widetilde p_{n,t,M} \|_2 ) \\
	& \leq  \frac 2 {NL}\sum_{t,n=1}^{L,N} ( \|\widetilde \varphi_{n,t-1,M} \|_2 + \|\widetilde p_{n,t,M} \|_2 ) \leq  \frac 2 {NL}\sum_{t,n=1}^{L,N} \sum_{k=1}^K ( | \widetilde \varphi_{n,t-1,M}(k)| + | \widetilde p_{n,t,M}(k)|), 
\end{align*}
where the last inequality follow from \eqref{eq:phi_p_tilde}. 
Applying \eqref{eq:phi_p_concentr}, we obtain 
\begin{align*}
& \prob{ \|\bb_M^e -\bb_\infty^e\|_{F} >  s}  \leq   \prob{ \frac 2 {NL}\sum_{t,n=1}^{L,N} \sum_{k=1}^K ( | \widetilde \varphi_{n,t-1,M}(k)| + | \widetilde p_{n,t,M}(k)|) >s }  \\
& \leq \sum_{t,n=1}^{L,N} \sum_{k} \left( \prob{ | \widetilde \varphi_{n,t-1,M}(k)| >\frac{s}{4K}  } + \prob{ | \widetilde p_{n,t,M}(k)| >\frac{s}{4K}  } \right) \leq  4NLK\exp\left( -\frac{Ms^2}{32 K^2+\frac{16sK}{3}} \right). 
\end{align*}
This completes the proof. 
\end{proof}

\begin{proof}[Proof of Theorem \ref{thm:concentrationT}] The proof is based on the Bernstein concentration inequalities for the normal matrices and the normal vectors. 

\textbf{Part (a).} Note that 
\begin{align*}
\|\widehat \bT_M - \bT \|_{F} & = \|\bA_M^\dagger \bb_M - \bA_\infty^{-1} \bb_\infty \|_{F} \leq  \|\bA_M^\dagger \bb_M - \bA_M^\dagger \bb_\infty\|_{F} +\| \bA_M^{\dag} \bb_\infty -    \bA_\infty^{-1} \bb_\infty \|_{F} \\
& \leq  \|\bA_M^\dagger \|_{op} \| \bb_M - \bb_\infty\|_{F} +\| \bA_M^{\dag} \|_{op}  \| ( \bA_\infty - \bA_M) \|_{op}\|\bT \|_{F}, 
\end{align*}
where the last inequality uses the fact that $\| ( \bA_M^{\dag}  -    \bA_\infty^{-1}) \bb_\infty \|_{F} = \| \bA_M^{\dag} (  \bA_\infty - \bA_M)  \bA_\infty^{-1} \bb_\infty \|_{F}  \leq \| \bA_M^{\dag} \|_{op}  \| ( \bA_\infty - \bA_M) \bT \|_{F}$. 
Hence, we have 
\[
\|\widehat \bT_M - \bT \|_F\leq \frac{1}{2} \epsilon + \frac{1}{2} \epsilon = \epsilon, \text{  on } E_1\cap E_2\cap E_3,  
\]
where we denote by $E_1,E_2,E_3$ the following events: 
\begin{align*}
 E_1 &: =\big \{\|\bA_M^\dag\|_{op} \leq C \big \} , \quad  \text{ with } C= 2 \|\bA_\infty^{-1}\|_{op} \\
 E_2 &: = \big \{ \| (\bA_M - \bA_\infty)\|_{op}\leq \frac{\epsilon}{2C\|\bT\|_F}  \big \} , \quad \quad 
 E_3 : = \big  \{ \| (\bb_M - \bb_\infty)\|_F \leq \frac{\epsilon}{2C}  \big \} .  
\end{align*}
Thus,  $\{\|\widehat \bT_M - \bT \|_F> \epsilon \}\subset E_1^c\cup E_2^c\cup E_3^c $.  Then, if we can prove the following bounds 
\begin{equation}\label{eq:E123}
\prob{E_1^c} <\frac{\delta}{3}, \quad  \prob{E_2^c} <\frac{\delta}{3}, \quad\prob{E_3^c} <\frac{\delta}{3}
\end{equation}
for $M\geq M_{\epsilon, \delta}$, we can conclude \eqref{eq:non-aysmp} by noting that 
\begin{align*}
\P\{\|\widehat \bT_M - \bT \|_F>  \epsilon \} < \P\{ E_1^c\cup E_2^c\cup E_3^c\} \leq \P\{ E_1^c\} + \P\{ E_2^c\} +  \P\{ E_3^c\} <\delta. 
\end{align*}

In the following, we prove the three bounds in \eqref{eq:E123} by Bernstein's inequalities. Note that 
$\|\bA_M^\dag\|_{op}=  \lambda_{min}(\bA_M)^{-1}$ and $\|\bA_\infty\|_{op} = \lambda_{min}(\bA_\infty)$. 
Thus, 
\begin{align*}
E_1^c & = \big \{\lambda_{min}(\bA_M)^{-1} > 2 \lambda_{min}(\bA_\infty)^{-1} \big \}  
   =  \{ \lambda_{min}(\bA_M) < \frac{1}{2}  \lambda_{min}(\bA_\infty)\} \\
 & \subset \big \{ |\lambda_{min}(\bA_\infty) - \lambda_{min}(\bA_M) | > \frac{1}{2}  \lambda_{min}(\bA_\infty)  \big \}.
\end{align*} Meanwhile, by Weyl's inequality, $|\lambda_{min}(\bA_\infty) - \lambda_{min}(\bA_M) |\leq \|\bA_\infty- \bA_M\|_{op}$. Hence,  by matrix Bernstein's inequality, 
\[
\P\{E_1^c\} <  \prob{ \|\bA_\infty- \bA_M\|_{op} >  \frac{1}{2}  \lambda_{min}(\bA_\infty)  } , 
\]
Similarly, 
$
\P\{E_2^c\} =  \prob{ \|\bA_\infty- \bA_M\|_{op} >  
\frac{\epsilon}{4\|\bT\|_F} \lambda_{min}(\bA_\infty) } $. 
Thus, with $s=\frac{ \lambda_{min}(\bA_\infty)}{2}  \min\big\{ 1, \frac{\epsilon}{2\|\bT\|_F}\big \} $, Lemma \ref{eq:concentrationAb} implies, 
\[
\max\big\{ \P\{E_1^c\}, \P\{E_2^c\}  \big\} <   \prob{ \|\bA_\infty- \bA_M\|_{op} >  s} < 2K \exp\left(-\frac{Ms^2/2}{1+s/3} \right). 
\]
Similarly, with $\alpha=\frac{\epsilon}{2C}  = \frac{\epsilon}{4}\lambda_{min}(\bA_\infty) $,  Lemma \ref{eq:concentrationAb} implies, 
\[
 \P\{E_3^c\} = \prob{ \|\bb_\infty- \bb_M\|_{F} >  \alpha} < 2K^2 \exp\left(-\frac{M\alpha^2}{8 K^2+4K\alpha/3} \right). 
\]
Hence, to obtain \eqref{eq:E123}, we set $M$ to satisfy both $2K \exp\left(-\frac{Ms^2/2}{1+s/3} \right)< \frac{\delta}{3}$ and $2K^2 \exp\left(-\frac{M\alpha^2}{8K^2+4K\alpha/3} \right)< \frac{\delta}{3}$, which lead to the lower bound for $M$.

\textbf{Part (b)}. The proof is similar to the above for Part (a), and we omit it here.   
\end{proof}

\paragraph{Acknowledgement} 
EB was supported by NSF-2106556, RW by NSF-2308120, MM and FL by FA9550-21-1-0317 and FA9550-23-1-0445, and FL by NSF-1913243 and NSF-2238486. The authors thank the anonymous referees, the Associate Editor, and Editor Davy Paindaveine for their comments.  
MM and FL would like to thank Andy Feinberg for discussions in epigenetics that provided partial motivation for introducing these types of PCA models.

\bibliographystyle{abbrv}

\newcommand{\etalchar}[1]{$^{#1}$}

\end{document}